\documentclass[12pt,a4paper,reqno]{amsart}
\usepackage[latin1]{inputenc}
\usepackage{amsmath}
\usepackage{amsfonts}
\usepackage{amssymb}
\usepackage{amsmath,bbm,amssymb,amsxtra}
\usepackage{enumerate}
\usepackage{caption}
\allowdisplaybreaks
\usepackage{epsfig}
\usepackage{ae}
    \def\Xint#1{\mathchoice
    {\XXint\displaystyle\textstyle{#1}}%
    {\XXint\textstyle\scriptstyle{#1}}%
    {\XXint\scriptstyle\scriptscriptstyle{#1}}%
    {\XXint\scriptscriptstyle\scriptscriptstyle{#1}}%
    \!\int}

    \def\XXint#1#2#3{{\setbox0=\hbox{$#1{#2#3}{\int}$ }
    \vcenter{\hbox{$#2#3$ }}\kern-.6\wd0}}
    \def\dashint{\Xint-}

\theoremstyle{definition}

\newtheorem{lemma}{Lemma}[section]
\newtheorem*{lemma2}{Lemma}
\newtheorem{proposition}[lemma]{Proposition}

\newtheorem{theorem}[lemma]{Theorem}

\newtheorem{corollary}[lemma]{Corollary}

\newtheorem{remark}[lemma]{Remark}

\newtheorem{definition}[lemma]{Definition}

\newcommand{\prop}[1]{\begin{proposition}\label{#1}
\sl }
\newcommand{\eprop}{\end{proposition}}
\newcommand{\thm}[1]{\begin{theorem}\label{#1}
\sl }
\newcommand{\ethm}{\end{theorem}}
\newcommand{\cor}[1]{\begin{corollary}\label{#1}
\sl }
\newcommand{\ecor}{\end{corollary}}

\newcommand{\lem}[1]{\begin{lemma}\label{#1}
\sl }
\newcommand{\elem}{\end{lemma}}
\newcommand{\defin}[1]{\begin{definition}\label{#1}
\sl }
\newcommand{\edefin}{\end{definition}}
\newcommand{\beqno}{\begin{eqnarray*}}
\newcommand{\eeqno}{\end{eqnarray*}}
\newcommand{\beqla}[1] {\begin {eqnarray}\label{#1}}
\def\eeq {\end {eqnarray}}
\newcommand{\beq}{\begin {eqnarray}}
\newcommand{\real}{{\mathbb R}}

\newcommand{\integer}{{\mathbb Z}}
\newcommand{\R}{{\mathbb R}}
\newcommand{\nanu}{{\mathbb N}}

\newcommand{\complex}{{\mathbb C}}

\renewcommand{\:}{\colon}

\newcommand{\trace}{{\mathbb T\hskip-4pt{\mathbb R}} }
\newcommand{\integral}{{\mathbb I}}

\newcommand{\diam}{{\rm diam}\,}

\newcommand{\Lip}{{\rm Lip}\,}
\newcommand{\D}{\mathbb{D}}
\newcommand{\hlmax}{\mathcal{M}}
\newcommand{\bzero}{\bar 0}
\newcommand{\df}{\dot{\mathcal F}}
\newcommand{\dfp}{\dot{\mathcal F}^{s}_{p,q}(Z)}

\newcommand{\dm}{\dot{M}}
\newcommand{\dmp}{\dot{M}^s_{p,q}(Z)}
\newcommand{\gkx}{\Gamma_{\kappa,\xi}}
\newcommand{\gox}{\Gamma_{1,\xi}}
\newcommand{\J}{{\mathcal J}}
\newcommand{\jpk}{\J^{s,(\kappa)}_{p,q}(X)}
\newcommand{\jp}{\J^{s}_{p,q}(X)}
\newcommand{\jpe}{\J^{s}_{p,q}(E)}
\newcommand{\jzero}{\mathring{\J}^s_{p,q}(X)}
\newcommand{\differencespace}{\mathring{{\mathcal D}}\jp}

\DeclareMathOperator*{\esssup}{ess\,sup}

\title[Triebel-Lizorkin spaces]{Triebel-Lizorkin spaces on metric spaces via hyperbolic fillings}

\author[Bonk]{Mario Bonk}
\address{University of California, Los Angeles, Department of Mathematics,  Box 95155, Los Angeles, CA, 90095-1555, USA.}
\email{mbonk@math.ucla.edu}
\thanks{M.B.\ was partially supported by NSF grant  DMS-1162471.}

\author[Saksman]{Eero Saksman}
\address{Department of Mathematics and Statistics, University of Helsinki, PO~Box~68, FI-00014 Helsinki, Finland.}
\email{eero.saksman@helsinki.fi}
\thanks{E.S.\ was supported by the Finnish CoE in Analysis and Dynamics Research, and by the Academy of Finland, projects 113826 and 118765.}

\author[Soto]{Tom\'as Soto}
\address{Department of Mathematics and Statistics, University of Helsinki, PO~Box~68, FI-00014 Helsinki, Finland}
\email{tomas.soto@helsinki.fi}
\thanks{T.S.\ was supported by the Finnish CoE in Analysis and Dynamics Research, and by the V\"ais\"al\"a foundation.}


\date{November 21, 2014}


\begin{document}

\maketitle

\begin{abstract}
We give a new characterization of (homogeneous) Triebel-Lizorkin spaces $\dot{\mathcal F}^{s}_{p,q}(Z)$ in the smoothness range $0 < s < 1$ for a fairly general class of metric measure spaces $Z$. The characterization uses Gromov hyperbolic fillings of $Z$. This gives a short proof of the quasisymmetric invariance of these spaces in case $Z$ is $Q$-Ahlfors regular and $sp = Q > 1$. We also obtain first results on complex interpolation for these spaces in the framework of doubling metric measure spaces.
\end{abstract}

\section{Introduction}\label{se:introduction}

The homogeneous Triebel-Lizorkin spaces $\dot{F}^s_{p,q}(\real^d)$ form a natural scale of function spaces that contains many of the important classical spaces such as the standard Sobolev, Hardy-Sobolev and BMO-type spaces. In particular, one has $\dot{F}^0_{p,2}(\real^d) = L^p(\real^d)$ and $\dot{F}^1_{p,2}(\real^d) = \dot{W}^{1,p}(\real^d)$. In addition, the real variable Hardy spaces $H^p(\real^d)$, the Hardy-Sobolev spaces $H^{s,p}(\real^d)$ with $s>0$ and the space of bounded mean oscillation $BMO(\real^d)$ are also included in the Triebel-Lizorkin scale. More details can be found in \cite{T,FJ,G}, for example.

Recently, many aspects of analysis, including the theory of Sobolev spaces, have been successfully carried over to the setting of metric measure spaces with many applications (see, for example, \cite{H2,S,Ch}). Besov and Triebel-Lizorkin spaces on metric spaces have been studied in many works, notably in \cite{HS,HMY,KYZ}. In \cite{KYZ} Triebel-Lizorkin spaces were defined on general metric measure spaces in terms of so-called \emph{fractional Haj\l asz gradients}; alternative characterizations can be found in \cite{GKZ}.

In the present paper, we give a new characterization of Triebel-Lizorkin spaces on a doubling metric measure space $Z$ in terms of ``Poisson extensions'' of locally integrable functions; the role of the upper half-space is taken by a hyperbolic filling $X$ of $Z$ (see \cite{BP}, \cite{BuS}, or Section \ref{se:definitions} below). The definition is in some sense analogous to the treatment of Besov spaces in \cite{BP} (see also \cite{BS}). In this way, one obtains natural descriptions of the Triebel-Lizorkin space $\dfp$ both as a closed subspace of $L^1_{\text{loc}}(Z)$ as well as a quotient space of a certain sequence space intrinsically defined on $X$. These spaces coincide with the Triebel-Lizorkin spaces introduced in \cite{KYZ} for most indices in the smoothness range $0 < s < 1$. In particular, they generalize the standard Triebel-Lizorkin spaces in case $Z=\real^d$.

A main motivation for the new definition is its intrinsic conformally invariant nature which gives the quasi-conformal invariance of these spaces as an easy application. The new definition is also amenable to complex interpolation results, as will be discussed below. Moreover, this approach in conjunction with a suitable retraction result (see Proposition \ref{pr:retraction} below) leads to new and general trace theorems for Sobolev, Besov and Triebel-Lizorkin spaces in the setting of Ahlfors regular metric spaces; this question will be studied in a separate paper \cite{SS}. 

For the reader's convenience we first recall one of the classical definitions of the norm on $\dot{F}^s_{p,q}(\real)$ (with $0 < s < 1$ and $1 < p,\,q < \infty$, say). Consider a smooth function $f$ on $\real^d$ and let $P_k(f)$ for $k\in\integer$ stand for a suitable approximation of $f$ that yields a ``good'' approximation of $f$ at levels of oscillation up to scale $2^{-k}$. Then, by defining a level $2^{-k}$ approximation of the derivative $d_kf:=P_k(f)-P_{k-1}(f)$, we may write
\[
  f = \sum_{k=-\infty}^\infty d_kf
\]
(at least up to an additive constant), and the Triebel-Lizorkin norm of $f$ is obtained as
\[
  \|f\|_{\dot{F}^s_{p,q}(\real^d)} = \bigg(\int_{\real^d} \big\| \big(2^{ks}d_kf(x)\big)_{k\in\integer}\big\|^p_{\ell^q} \,dx \bigg)^{1/p}.
\]
Here $\|(y_k)\|_{\ell^q}:=\big(\sum_{k\in\integer}|y_k|^q\big)^{1/q}$ stands for the $\ell^q$-norm of a given sequence $(y_k)_{k\in\integer}.$

An example of a valid approximation is the standard Fourier-analytic definition; in this case one takes $\widehat{P_k(f)} = \varphi(2^{k}\cdot) \widehat{f}$, where $\varphi$ is a suitable function supported in the unit ball. Another suitable approximation is obtained by using the partial sums of a wavelet decomposition.

\newcommand{\dist}{{{\rm dist}}}

In order to explain the main idea of our definition for $\dfp$ in the setting of metric measure spaces $Z$, we still focus $Z=\real^d$ for a moment. Let $Q_0:=[0,1]^d$ be the unit cube in $\real^d$, and $Q_{j,k}= 2^{-k}(j+Q_0)$ for $k\in \integer$ and $j\in \integer^d$. Then $\{ Q_{j,k} : j\in \integer^d,\, k\in\integer \}$ is the collection of dyadic cubes in $\real^d$. We consider the upper half-space $\real^{d+1}_+=\real^d\times\real_+=\{(z,y):z\in \real^d,\, y>0\}$ and let $\widehat Q_{j,k}:=Q_{j,k}\times[2^{-k-1}, 2^{-k}]\subset\real^{d+1}_+$. Then $\{ \widehat Q_{j,k} : j\in \integer^d,\, k\in\integer \}$ is a Whitney-type decomposition of $\real^{d+1}_+$.

For $f \in L^1_{\rm{loc}}(\real^d)$ we can now define the approximation $Pf$ (a discrete analog of the Poisson extension) as a function on the set of centers of these Whitney cubes: if we denote by $y_Q$ the center of a Whiney cube $\widehat Q$ associated with a dyadic cube $Q\subset \real^d$, then we set
\[
  Pf(y_Q) := f_Q := \frac{1}{m_d(Q)}\int_Qf\, dm_d,
\]
where $m_d$ denotes Lebesgue measure on $\R^d$. Moreover, a discrete derivative of $Pf$ (or rather its absolute value) is obtained by setting
\beqla{eq:dP}
  |d(Pf)(y_Q)| := \sum_{\widehat Q' \cap \widehat Q\ne \emptyset}|f_{Q'}-f_Q| \quad \textrm{for} \quad Q \in \{ Q_{j,k} : j \in \integer^d,\,k\in\integer\}.
\eeq
The quantity $|d(Pf)(y_{Q_{j,k}})|$ can be thought as a version of $d_kf$ on the cube $Q_{j,k}$. 

In order to define the Triebel-Lizorkin norm, we still need an analog of the quantity $\| \big( 2^{ks}d_kf(x) \big)_{k\in\integer}\|_{\ell^q}$. To this end, for given $\xi \in \real^d$ we consider the standard non-tangential cone $\Gamma_{\xi}$ in the upper half-space with tip at $\xi$ defined as 
\[
  \Gamma_{\xi} = \{ (z,y)\in \real^d \times \real_+ \,:\, |z-\xi |<y\}.
\]
We also set $\ell(y_Q) := k$ if $Q \in \{ Q_{j,k} \,:\, j \in \integer^d\}$. A natural replacement for the quantity $\| \big(2^{ks}d_kf(\xi)\big)_{k\in\integer}\|_{\ell^q}$ is then given by
$
  \| \{ 2^{\ell(y_Q)s} |d(Pf)(y_Q)| : y_Q \in \Gamma_\xi \} \|_{\ell^q}
$, 
because it selects values of the discrete derivative $d(Pf)$ only from Whitney cubes that lie above $\xi$ (in particular, it selects a uniformly bounded number of values $2^{\ell(y_Q)} |d(Pf)(y_Q)|$ from each level $k$). With these identifications, we are led to the definition of a (semi-)norm
\beqla{eq:real_d}
  \|f\|_{\df^s_{p,q}(\real^d)} := \bigg(\int_{\real^d} \big\|\{ 2^{\ell(y_Q)} |d(Pf)(y_Q)| : y_Q \in \Gamma_\xi \}\big\|^p_{\ell^q} d\xi \bigg)^{1/p},
\eeq
and of $\df^s_{p,q}(\real^d)$ as the space of all functions $f\in L^1_{\rm loc}(\real^d)$ for which $\|f\|_{\df^s_{p,q}(\real^d)}<\infty$.

This definition now easily carries over to the setting of metric measure spaces $(Z,d,\mu)$: one replaces the dyadic cubes by a collection of balls $\{ B \}$ obtained from suitable covers of $Z$ on each scale $2^{-k}$, one uses the averages $f_B$ to define $Pf$ and obtains $d(Pf)$ as in \eqref{eq:dP}. The cone $\Gamma_\xi$ at $\xi \in Z$ is simply defined as the set of the centers of the balls $B$ that contain $\xi$. The definition \eqref{eq:real_d} then generalizes to the setting of arbitrary metric measure spaces $(Z,d,\mu)$.

In the case of the Euclidean space $\real^d$, this definition is equivalent to the definition of the Triebel-Lizorkin spaces $\dot{F}^s_{p,q}(\real^d)$ for $0 < s < 1$, $d/(d+s) < p < \infty$, and $d/(d+s) < q \leq \infty$. For most indices the definition is also equivalent to the one given in terms of fractional Haj\l asz gradients in \cite{KYZ}: roughly speaking, one looks at the mixed $L^p(\ell^q)$ norms of sequences $(g_k)_{k \in \integer}$ of measurable functions $g_k \: Z \to [0,\infty]$ satisfying
\[
  |f(\xi) - f(\eta)| \leq d(\xi,\eta)^s \big(g_k(\xi) + g_k(\eta)\big), 
\]
whenever $\xi$, $\eta \in Z$ and $2^{-k-1} \leq d(\xi,\eta) < 2^{-k}$. 

Our definition is slightly harder to state and remember than the one using Haj\l asz gradients, but it has the advantage that it uses a linear discrete derivative operator in order to define the norm directly, whereas the Haj\l asz gradients are in general not uniquely determined. This feature, together with the fact that our approach yields an identification of the spaces $\dfp$ with sequence spaces, also allows one to extend the standard result on complex interpolation between Triebel-Lizorkin spaces to cover the case of fairly general metric spaces $Z$. More precisely, in Theorem \ref{th:interpolation-f} we establish the interpolation formula
\[
  \big[\df^{s_0}_{p_0,q_0}(Z),\df^{s_1}_{p_1,q_1}(Z)\big]_{\theta} = \dfp
\]
for a wide range of indices, where $s$, $p$ and $q$ on the right-hand side are related to the indices on the left-hand side in the expected way. Among other new results in the general metric setting, we generalize classical embedding results for these spaces and verify the density of Lipschitz functions, see Proposition \ref{pr:jawerth-embedding} and Theorem \ref{th:approximations} below.

The quantity $d(Pf)$ is a natural analog of a ``hyperbolic derivative'' of the Poisson extension, and one may consider the centers of the balls in the collection $\{B\}$ associated with $Z$ as vertices of a naturally defined Gromov hyperbolic graph $X$, called the hyperbolic filling of $Z$. In this setup, it is known that quasisymmetric self-maps of $Z$ correspond to quasi-isometries of the graph $X$. Hence our definition for the spaces $\dfp$ is ideally set up to study the quasisymmetric invariance properties of these spaces; indeed, we obtain a concise approach to the important quasiconformal invariance result of \cite{KYZ}, whose original proof is a technical ``tour de force''. The new proof of this result is included in Section \ref{se:invariance}.

We note here that an analogous definition using hyperbolic fillings can be given  for Besov spaces, and most of the results of the present paper (outside Section \ref{se:invariance}) hold true for these spaces as well; see \cite{So}.

The structure of the paper is as follows: Section \ref{se:definitions} contains the description of the hyperbolic filling $X$ of a doubling metric measure space $Z$ and definitions of the Poisson extension $Pf \: X \to \complex$ of a given function $f \in L^1_{\rm loc}(Z)$ as well as the discrete derivative $d(Pf)$ of the extension. The definition of $\dfp$ is then obtained in terms of $d(Pf)$ in a manner analogous to the one in classical harmonic analysis in case $Z = \real^d$. To be slightly more precise, the space $\dfp$ is defined by the quasi-norm
\[
  f \mapsto \|d(Pf)\|_{\jp},
\]
where $\jp$ is the sequence space alluded to above; see Definition \ref{de:fspace} for details. The section also contains some important auxiliary results, especially Theorem \ref{th:trace-properties} on the commutativity of the trace operator and the Poisson extension.

In Section \ref{se:equivalence} we prove that in the smoothness range $0 < s <1$, our definition of the Triebel-Lizorkin spaces is equivalent to the one in \cite{KYZ}; in particular, this gives the equivalence of our definition  with the usual Fourier analytical definition in case $Z = \real^d$. In addition, we prove that Lipschitz functions are dense in $\dfp$ when $q < \infty$.

For simplicity, we assume in Sections \ref{se:definitions} and \ref{se:equivalence} that the metric space $Z$ is bounded. In Section \ref{se:unbounded} we treat the case of unbounded $Z$. The definition of the spaces remains unchanged and here we indicate the minor modifications needed to generalize the results of Sections \ref{se:definitions} and \ref{se:equivalence} to this situation. Section \ref{se:invariance} contains the short proof of the quasisymmetric invariance result mentioned above, and finally our new results on complex interpolation are given in Section \ref{se:interpolation}.

Let us end this section by introducing some notation. On any metric space $(Z,d)$, we write $B_d(\xi,r):=\{\eta \in Z : d(\eta,\xi) < r\}$ for the open ball centered at $\xi \in Z$ with radius $r > 0$, or just $B(\xi,r)$ if the metric is obvious from the context. For an arbitrary ball $B$ with a distinguished center point $\xi$ and radius $r > 0$, we let $\lambda B := B(\xi,\lambda r)$ for $\lambda > 0$. If $\mu$ is a measure on $Z$, $E$ is a subset of $Z$ and $f$ is a complex-valued function on $Z$, we write
\[
  \dashint_E f d\mu := \frac{1}{\mu(E)} \int_E f d\mu,
\]
whenever the latter quantity is well-defined. Instead of the usual definition, we shall denote by $L^1_{\rm loc}(Z)$ the vector space of complex-valued measurable functions on $Z$ that are integrable on \emph{bounded measurable subsets of $Z$}. For any two non-negative functions $f$ and $g$ with the same domain, the notation $f \lesssim g$ means that $f \leq C g$ for some finite constant $C>0$, independent of certain parameters that will be obvious from the context. The notation $f \approx g$ means that $f \lesssim g$ and $g \lesssim f$. We will write $a\land b$ for the minimum of two real numbers $a$ and $b$.

\section{Definitions and basic results}\label{se:definitions}

Let $(Z,d)$ be a metric space, and let $\mu$ be a Borel regular measure on $Z$ such that the $\mu$-measure of every open ball is positive and finite. We assume that the metric measure space $(Z,d,\mu)$ satisfies the following \emph{doubling property}: $\mu(2B)\leq c\mu(B)$ for all balls $B$, where the finite constant $c$ is independent of $B$. It is an elementary consequence of this property that
\beqla{eq:doubling}
  \mu\big(B(\xi,\lambda r) \big) \leq C \lambda^Q \mu\big(B(\xi,r)\big)
\eeq
for some constants $Q > 0$ and $C \geq 1$ whenever $\xi \in Z$, $r > 0$, and $\lambda \geq 1$. The number $Q$ will be fixed from now on.

In this section, we will consider spaces with finite diameter. This is done for the sake of simplicity, and the simple modifications needed in the unbounded case are considered later in Section \ref{se:unbounded}. We may also assume that $\diam Z = 1$ for convenience.

We denote by $X$ the ``hyperbolic filling'' with respect to $(Z,d,\mu)$ as defined in \cite{BP} and \cite{BS} (see also \cite[Chapter 6]{BuS}). More precisely, for any integer $n \geq 1$, we choose a maximal set of points $\{\xi_x\}_{x \in X_n}$ of $Z$, where $X_n$ is some index set, so that $d(\xi_x,\xi_{x'}) \geq 2^{-n-1}$ whenever $x$, $x' \in X_n$ and $x \neq x'$. By the doubling property, the balls $B(x) := B(\xi_x,2^{-n})$ then have bounded overlap (uniformly in $n$), and it is easily seen that the balls $2^{-1} B(x)$ cover $Z$ and that the balls $4^{-1} B(x)$ are pairwise disjoint. We further let $X_0$ be a singleton set, whose element we shall denote by $\bzero$, and define $B(\bzero) := Z$. Write $|x| := n$ for all $x \in X_n$, $n \geq 0$. We then consider the disjoint union $X := \sqcup_{n \geq 0} X_n$, and denote by $(X,E)$ the graph such that the vertices $x$, $x' \in X$ are joined by an edge in $E$ if and only if $x\ne x'$, $||x|-|x'|| \leq 1$, and $B(x)\cap
 B(x') \neq \emptyset$.

\begin{figure}[hbt]\caption*{\textsc{Figure.} Illustration of the levels $X_{-1}$ through $X_3$ of a hyperbolic filling of the unit interval $[0,1]$.} \bigskip\bigskip \centering \includegraphics[scale=0.45]{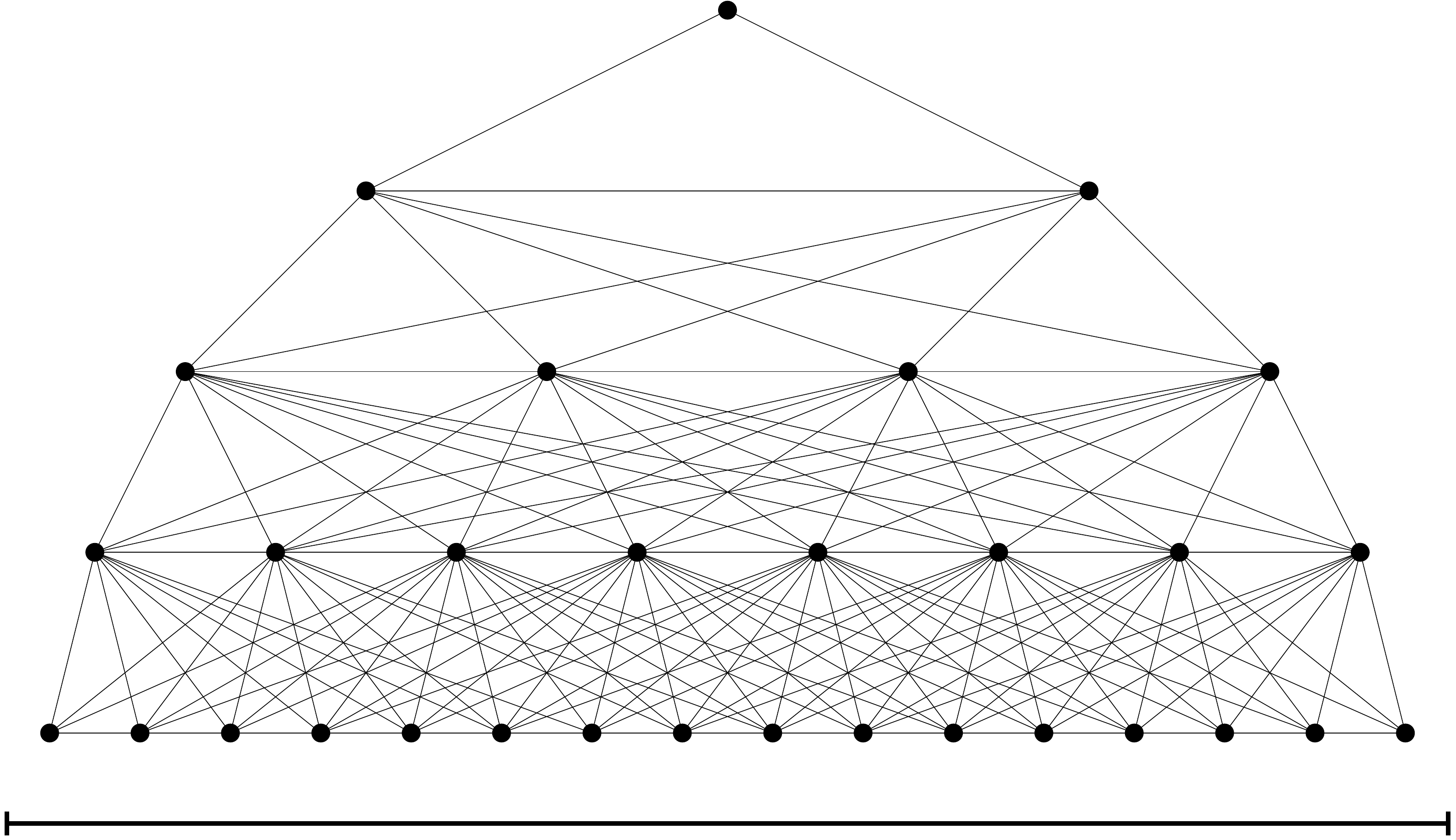} \end{figure}

If $x'\not=x$ is joined to $x$ by an edge, we say that $x'$ is a neighbor of $x$ and write $x'\sim x$. Under some mild additional assumptions on $Z$, the natural path metric on the graph $(X,E)$ makes it hyperbolic in the sense of Gromov so that its boundary at infinity coincides with $Z$ (see \cite{BP}); this is the reason why $X$ is called hyperbolic filling of $Z$. 

\newcommand{\degree}{\Delta(X)}
\newcommand{\degr}{{\rm deg}}

For a complex-valued function $u$ on the hyperbolic filling $X$, the \emph{discrete derivative} $du$ can either be defined%
\footnote{The definition of $du$ as a function on edges is as follows: if $(X,E)$ is defined as a directed graph, the value of $du$ on the dirrected edge from $x$ to $x'$ is given by $u(x')-u(x)$.}
on the set of all edges $E$ or, in an essentially equivalent way, as a vector-valued function on $X$ itself (taking values in a finite-dimensional vector space). We shall adopt the latter point of view and set
\[
  du(x):=\big(u(y_1)-u(x),\ldots, u(y_{\degr(x)})-u(x),0,\ldots\big)\in \complex^{\degree}
\]
for $x\in X$, where $y_1,\ldots ,y_{\degr(x)}\in X$ are the neighbors of $x$ listed in some (fixed) order, and $\degree:=\sup_{x\in X}m(x)<\infty$ stands for the maximal vertex degree of $(X,E)$; note that $\degree $ is finite as follows from the fact that $Z$ is doubling. In our definition the vector $du(x)$ is obtained by augmenting the vector $(u(y_1)-u(x),\ldots, u(y_{\degr(x)})-u(x))$ by zeroes, if needed, so that it has precisely $\degree$ coordinates. Thus the map $u\mapsto du$ is linear from the space of functions on $X$ to the space $(\complex^{\degree})^X$. Actually we will not make much use of the full derivative $du(x)$, but will just use its pointwise norm given by
\[
  |du|(x):=\big(\sum_{y\sim x}|u(y)-u(x)|^2\big)^{1/2}.
\]

The \emph{Poisson extension} $Pf \: X \to \complex$ of a function $f \in L^1_{\rm loc}(Z)$ is defined by
\[
  Pf(x) := \dashint_{B(x)} f \,d\mu
\]
for $x \in X$. 

The \emph{cone} $\gkx$ of width $\kappa \in [1,\infty)$ at a point $\xi \in Z$ is defined as the set
$
 \gkx := \left\{ x\in X : \xi \in \kappa B(x) \right\}.
$
One can view $\gkx$ as a discrete analog of a non-tangential approach region (or a Stoltz angle) at the point $\xi$ in the Euclidean case. 

For any set $A$ of complex numbers, we shall write $\|A\|_{\ell^q}$ for $(\sum_{a \in A} |a|^q)^{1/q}$ when $0 < q < \infty$ and for $\sup_{a \in A} |a|$ when $q = \infty$. Note that with $A$ fixed, $\|A\|_{\ell^q}$ is a non-increasing function of $q$.

We next introduce a family of sequence spaces that will subsequently be used to define the Triebel-Lizorkin spaces.

\defin{de:sspace}
For $s \in (0,\infty)$, $p \in (0,\infty)$, and $q \in (0,\infty]$, we define the sequence space $\jp$ as the quasi-normed space of functions $u \: X \to \complex$ such that
\[
  \|u\|_{\jp} := \Big( \int_{Z} \big\| \{ 2^{|x|s} |u(x)| : x \in \gox \} \big\|_{\ell^q}^p \,d\mu(\xi)\Big)^{1/p}<\infty.
\]
\edefin

One easily checks that the integrand in the above definition is measurable. By the monotonicity of $\ell^q$-quasi-norms, we have that
\[
  \J^{s}_{p,q'}(X) \subset \jp
\]
with a continuous embedding whenever $q' \leq q$. When $p = q$, a simple calculation yields
\beqla{eq:j-diagonal}
  \|u\|_{\J^{s}_{p,p}(X)} = \big\|\{ 2^{|x|s}\mu(B(x))^{1/p} |u(x)| : x\in X\}\big\|_{\ell^p}.
\eeq

Here is the first basic result on the structure of these spaces.

\prop{pr:quasi-banach}
\textnormal{(i)} $\jp$ is a quasi-Banach space for all admissible parameters. When $1 < p, q < \infty$, it is a reflexive Banach space.

\smallskip
\textnormal{(ii)} Whenever $0 < s < \infty$, $0 < p < \infty$ and $0 < q \leq \infty$, the quasi-norm
\[
  \|u\|_{\jpk} := \Big( \int_{Z} \big\| \{ 2^{|x|s} |u(x)| : x \in \gkx \} \big\|_{\ell^q}^p \,d\mu(\xi)\Big)^{1/p}
\]
is equivalent to the quasi-norm of $\jp$ for all $\kappa \in [1,\infty)$.
\eprop

Before the proof, let us formulate the following auxiliary result, which will play an important role throughout the paper. For the proof, we refer to \cite[Theorem 1.2]{GLY}.

\begin{lemma2}[Fefferman-Stein maximal theorem for doubling metric measure spaces]\sl
Let $p \in (1,\infty)$ and $q \in (1,\infty]$. For any sequence $\vec{f} := \{f_k\}_{k \in \integer}$ of measurable functions $f_k \: Z \to \complex$, we have
\[
  \bigg(\int_Z \Big( \sum_{k \in \integer} \hlmax f_k (\xi)^q \Big)^{p/q} \,d\mu(\xi) \bigg)^{1/p} \lesssim \bigg(\int_Z \Big( \sum_{k \in \integer} |f_k (\xi)|^q \Big)^{p/q} \,d\mu(\xi) \bigg)^{1/p},
\]
(with an obvious modification for $q = \infty$), 
where $\hlmax$ stands for the Hardy-Littlewood maximal operator on $(Z,d,\mu)$ and the implicit constant is independent of $\vec{f}$.
\end{lemma2}

\begin{proof}[Proof of Proposition \ref{pr:quasi-banach}]
(i) Let $r := \min(1,p,q)$. To show that a Cauchy sequence $(u_k)_{k \geq 1}$ in $\jp$ converges in $\jp$, we can assume that $\|u_{k+1} - u_k\|_{\jp} \leq 2^{-k}$ for all $k$. We have
\[
  |u_{k+1}(x) - u_k(x)| \lesssim 2^{-|x|s}\mu\big(B(x)\big)^{-1/p} \|u_{k+1} - u_k\|_{\jp}
\]
for all $x$ and $k$, and hence the limit
\[
  u(x) := u_1(x) + \sum_{k=1}^{\infty}\big( u_{k+1}(x) - u_k(x) \big)
\]
exists pointwise in $X$. In particular, we have
\[
  |u(x) - u_{N}(x)| \leq \sum_{k\geq N}|u_{k+1}(x) - u_{k}(x)|
\]
for all $x \in X$ and $N \in \nanu$. Since the function $v \mapsto \|v\|_{\jp}^r$ is clearly subadditive, this gives 
\[
  \|u-u_N\|_{\jp} \leq \Big( \sum_{k \geq N} \|u_{k+1} - u_k\|_{\jp}^r \Big)^{1/r}\to 0 \quad \text{as $N\to \infty$}.
\]

Now if $1 < p,q < \infty$, the vector-valued $L^p$-space $L^p(\ell^q(X,w_s))$ over $(Z,\mu)$, where $w_s$ stands for the weight $x \mapsto 2^{|x|s}$ on $X$, is known to be a reflexive Banach space. By definition, $u \mapsto ( u(x) \chi_{ B(x)})_{x\in X}$ is an isometric isomorphism of $\jp$ onto a subspace of $L^p(\ell^q(X,w_s) )$. By the first part of this proof, this subspace is closed, and so $\jp$ is reflexive.

(ii) It suffices to show that
\[
  \|u\|_{\jpk} \lesssim \|u\|_{\jp}
\]
for all sequences $u$ on $X$, with an implicit constant independent of $u$.

For arbitrary $u$, define the functions $U_k \: Z \to [0,\infty)$, $k \in \nanu_0$, by
\[
  U_k(\xi) := \big\| \{ 2^{|x|s} |u(x)| : x \in \gox\cap X_k \} \big\|_{\ell^q},
\]
so that
\[
  \|u\|_{\J^{s,(1)}_{p,q}(X)} = \bigg( \int_Z \| \{U_k(\xi) : k \geq 0\} \|_{\ell^q}^p \,d\mu(\xi)\bigg)^{1/p}.
\]
Now choose $r > 0$ so that $r < \min(p,q)$. If $\xi \in Z$ and $x \in \gkx \cap X_k$, we have
\[
  2^{|x|s}|u(x)| \leq \Big( \dashint_{B(x)} U_k(\eta)^r \,d\mu(\eta) \Big)^{1/r} \lesssim \Big( \dashint_{\kappa B(x)} U_k(\eta)^r \,d\mu(\eta) \Big)^{1/r} \lesssim \hlmax\big(U^r_k\big)(\xi)^{1/r}.
\]
Since $\#(\gkx \cap X_k)$ is bounded uniformly in $\xi$ and $k$, we thus have
\[
  \big\| \big\{2^{|x|s} |u(x)| : x \in \gkx \big\}\big\|_{\ell^q} \lesssim \big\| \big\{ \hlmax\big(U^r_k\big)(\xi) : k \geq 0 \big\}\big\|_{\ell^{q/r}}^{1/r}.
\]
Since $q/r$, $p/r > 1$, the Fefferman-Stein maximal theorem gives
\begin{align}
  \|u\|_{\jpk}^p\;\notag 
& \lesssim \;\int_Z \big\| \big\{ \hlmax\big(U^r_k\big)(\xi) : k \geq 0 \big\}\big\|_{\ell^{q/r}}^{p/r} \,d\mu(\xi) \notag \\
& \lesssim \;\int_Z \big\| \big\{ U_k(\xi)^r : k \geq 0 \big\}\big\|_{\ell^{q/r}}^{p/r} \,d\mu(\xi) \label{eq:fefferman-stein} \\
& = \;\|u\|_{\jp}^p. \qedhere
\end{align}
\end{proof}

We shall now examine the boundary behavior of functions in certain sequence spaces on $X$. We will first formulate an analog to \cite[Lemma 4.1]{BS}. To this end, fix a collection $(\psi_x)_{x\in X}$ of non-negative Lipschitz functions on $X$ such that $\psi_x$ is supported on $B(x)$ for all $x$, $(\psi_x)_{x\in X_n}$ is a partition of unity for all $n \geq 0$, and $\Lip \psi_x \lesssim 2^{|x|}$ for all $x$. The existence of such partitions is easy to verify. For a complex-valued function $u$ on $X$ we then define 
\[
  T_n u := \sum_{x \in X_n} u(x)\psi_x \quad \text{for $n \in \nanu_0$.}
\]

\lem{le:trace}
Suppose that $0 < s < \infty$, $Q/(Q+s) < p < \infty$, and $0 < q \leq \infty$. If $u$ is a sequence on $X$ such that $|du| \in \jp$, then the limit (the \textnormal{trace function} of $u$ on the boundary $Z$)
\[
  \trace u := \lim_{n\to\infty} T_n u
\]
exists in $L^{1}(Z)$ and pointwise $\mu$-almost everywhere, and
\[
  \|\trace u - u(\bzero)\|_{L^1(Z)} \lesssim \|du\|_{\jp}.
\]
\elem

The operator $\trace$ depends on the exact choice of the family $(\psi_x)_{x \in X}$. Since the choice is fixed in the sequel, we have not taken this dependence into account in the notation $\trace$.

\begin{proof}
Let $Q/(Q+s) < r < \min(p,1)$ and take $\epsilon \in (0,s)$ so that $r = Q/(Q+\epsilon)$. Now, if $x \in X_n$ for some $n$ and $\xi \in B(x)$, we have
\[
  |T_{n+1}u(\xi) - T_{n}u(\xi)| = \Bigg|\sum_{\substack{y'\in X_{n+1}\\\xi \in B(y')}}\Big(u(y') - u(x)\Big)\psi_{y'}(\xi) - \sum_{\substack{y\in X_{n}\\\xi \in B(y)}}\Big(u(y) - u(x)\Big)\psi_{y}(\xi)\Bigg|,
\]
where the latter quantity can be estimated from above by a constant times $|du(x)|$. This follows from the bounded overlap of the balls corresponding to the vertices at a fixed level of $X$. We conclude that 
\[
  \int_{Z} |T_{n+1}u(\xi) - T_{n}u(\xi)|\, d\mu(\xi) \lesssim \sum_{x\in X_n} \mu\big(B(x)\big)|du(x)|.
\]
Summing up over $n \geq 0$ and using the fact that $r < 1$, we have
\[
  \sum_{n \geq 0} \int_{Z} |T_{n+1}u(\xi) - T_{n}u(\xi)| \,d\mu(\xi) \lesssim \bigg( \sum_{x\in X} \Big[\mu\big(B(x)\big)|du(x)|\Big]^r \bigg)^{1/r}.
\]
Since $r = 1 - (\epsilon/Q)r$, the latter quantity is equal to 
\[
  \bigg( \sum_{x\in X} \Big[\mu\big(B(x)\big)^{1/r - (\epsilon/Q)}|du(x)|\Big]^r \bigg)^{1/r}.
\]
By the doubling property, we have $\mu(B(x))^{-(\epsilon/Q)} \lesssim 2^{|x|\epsilon} \mu(B(\bzero))^{-(\epsilon/Q)} \approx 2^{|x|\epsilon}$. So by \eqref{eq:j-diagonal}, the quantity above is controlled by a constant times $\|du\|_{\J^{\epsilon}_{r,r}(X)}$. Since $s > \epsilon$ and $\#\big(\gox \cap X_n\big) \approx 1$ uniformly in $\xi$ and $n$, we have
\[
  \big\| \{ 2^{|x|\epsilon} |du(x)| : x \in \gox \} \big\|_{\ell^r} \lesssim \big\| \{ 2^{|x|s} |du(x)| : x \in \gox \} \big\|_{\ell^q}
\]
for all $\xi \in Z$ regardless of the relationship between $r$ and $q$, so that
\[
  \sum_{n \geq 0} \int_{Z} |T_{n+1}u(\xi) - T_{n}u(\xi)| \,d\mu(\xi) \lesssim \|du\|_{\J^{s}_{r,q}(X)},
\]
and since $p > r$ and $\mu(Z) < \infty$, we can finally use H\"older's inequality to estimate the latter quantity by a constant times $\|du\|_{\jp}$. Hence
\[
  \sum_{n \geq 0} \int_{Z} |T_{n+1}u(\xi) - T_{n}u(\xi)| \,d\mu(\xi) \lesssim \|du\|_{\jp} < \infty.
\]
Both claims then follow from this and the fact that $T_0 u \equiv u(\bzero)$.
\end{proof}

Theorem \ref{th:trace-properties} below is the counterpart to \cite[Lemma 4.4]{BS}. Before stating it, let us formulate the following basic auxiliary result, which will be applied several times throughout this paper.

\prop{pr:basic-operator}
Let $0 < s < \infty$, $Q/(Q+s) < p < \infty$, and $Q/(Q+s) < q \leq \infty$. Then the operator $T$ defined by
\[
  Tu(x) := \sum_{\substack{|y| \geq |x|\\ B(y)\cap B(x) \neq \emptyset}} \frac{\mu\big(B(y)\big)}{\mu\big(B(x)\big)} u(y)
\]
is well-defined and bounded on $\jp$.

The conclusion continues to hold for the operator $u\mapsto (Tu)\circ \Psi$ whenever $\Psi \: X \to X$ is a mapping with $B(\Psi(x)) \cap B(x) \neq \emptyset$ for all $x$ and such that there exists  $\sigma \geq 0$ so that $|x| - \sigma \leq |\Psi(x)| \leq |x| + \sigma$ for all $x$.
\eprop

\begin{proof}
Let $u \in \jp$ and write
\[
  U_k(\xi) := \big\| \{ 2^{|x|s} |u(x)| : x \in \gox \cap X_k\} \big\|_{\ell^q},
\]
so that
\[
  \|u\|_{\jp} = \bigg( \int_Z \| \{ U_k(\xi) : k \geq 0\}\|_{\ell^q}^p \,d\mu(\xi) \bigg)^{1/p}.
\]

For $\xi \in Z$ and $x \in \gox$, we then have
\begin{align}
  2^{|x|s} T(|u|)(x)\;
& = \;2^{|x|s} \sum_{\substack{|y| \geq |x|\\ B(y)\cap B(x) \neq \emptyset}} \frac{\mu\big(B(y)\big)}{\mu\big(B(x)\big)} |u(y)| \notag\\
& = \;\sum_{k \geq |x|} 2^{(|x| - k)s} \sum_{\substack{y \in X_k \\ B(y) \cap B(x) \neq \emptyset}} \frac{\mu\big(B(y)\big)}{\mu\big(B(x)\big)} 2^{|y| s}|u(y)| \notag\\
& =: \;\sum_{k \geq |x|} S_k(x) \label{eq:sk-1}.
\end{align}

If we choose $r$ so that $Q/(Q+s) < r < \min(1,p,q)$, the doubling property implies
\[
  \Big(\frac{\mu\big(B(y)\big)}{\mu\big(B(x)\big)}\Big)^{1-1/r} \lesssim 2^{(|x|-|y|)(Q - Q/r)}
\]
for all $y$ such that $|y| \geq |x|$ and $B(y) \cap B(x) \neq \emptyset$; so for all $k \geq |x|$ the quantity $S_k(x)$ can be estimated from above by a constant times
\[
  2^{(|x| - k)(Q+s - Q/r)} \sum_{\substack{y \in X_k \\ B(y) \cap B(x) \neq \emptyset}} \bigg(\frac{\mu\big(B(y)\big)}{\mu\big(B(x)\big)}\bigg)^{1/r} 2^{|y|s}|u(y)|.
\]
Using the fact that $r < 1$ together with the bounded overlap of the balls corresponding to the elements of $X$ at any fixed level, we further get
\begin{align}
  S_k(x) \;
& \lesssim \; 2^{(|x| - k)(Q+s - Q/r)} \bigg( \sum_{\substack{y \in X_k \\ B(y) \cap B(x) \neq \emptyset}} \frac{\mu\big(B(y)\big)}{\mu\big(B(x)\big)} \big[2^{|y|s}|u(y)|\big]^r \bigg)^{1/r} \notag\\
& \leq\; 2^{(|x| - k)(Q+s - Q/r)} \bigg( \frac{1}{\mu\big( B(x)\big)}\sum_{\substack{y \in X_k \\ B(y) \cap B(x) \neq \emptyset}} \int_{B(y)} U_k(\eta)^r \,d\mu(\eta) \bigg)^{1/r} \notag\\
& \lesssim \;2^{(|x| - k)(Q+s - Q/r)} \bigg( \frac{1}{\mu\big( B(x)\big)} \int_{3B(x)} U_k(\eta)^r\, d\mu(\eta) \bigg)^{1/r} \notag\\
& \lesssim \;2^{(|x| - k)(Q+s - Q/r)} \hlmax \big( U_k ^r \big)(\xi)^{1/r}. \label{eq:sk-2}
\end{align}

As $Q+s-Q/r > 0$, using the subadditivity of $t \mapsto t^q$ if $q \leq 1$ or H\"older's inequality otherwise, we thus have
\begin{align*}
  \sum_{x \in \gox} \big[2^{|x|s} T(|u|)(x)\big]^q \;
& \lesssim \;\sum_{x \in \gox} \sum_{k \geq |x|}2^{(|x|-k)(Q+s-Q/r)(q\land 1)} \hlmax \big( U_k ^r \big)(\xi)^{q/r} \\
& \approx \;\sum_{k \geq 0} \hlmax \big( U_k ^r \big)(\xi)^{q/r}
\end{align*}
(with an obvious modification for $q = \infty$). Since $r < \min(p,q)$, one can obtain the first claim using the Fefferman-Stein maximal theorem as in \eqref{eq:fefferman-stein}.

The second claim readily follows from the first one and the fact that the composition operator $v \mapsto v\circ \Psi$ is bounded on $\jp$, which is easily seen using the fact that $X$ has bounded valency.
\end{proof}

\thm{th:trace-properties}
Let $0 < s < \infty$, $Q/(Q+s) < p < \infty$, and $Q/(Q+s) < q \leq \infty$.

\smallskip 
\textnormal{(i)} If $u \in \jp$, then $\trace u = 0$ $\mu$-almost everywhere.

\smallskip 
\textnormal{(ii)} If $u$ is a sequence on $X$ such that $|du| \in \jp$, then $u - P\trace u \in \jp$ and
\[
  \|u-P\trace u\|_{\jp} \lesssim \|du\|_{\jp}.
\]
In particular, $\trace u = \trace (P\trace u)$ pointwise $\mu$-almost everywhere.

\smallskip 
\textnormal{(iii)} If $f\in L^1(Z)$ and $|d(Pf)| \in \jp$, then $\trace (Pf) = f$ (with convergence in $L^1(Z)$ and pontwise $\mu$-almost everywhere).
\ethm

\begin{proof}
(i) For almost all $\xi \in Z$, the quantity $c_\xi := \sup_{x \in \gox} 2^{|x|s} |u(x)|$ is finite. For these $\xi$ we have $|T_n u(\xi)| \lesssim c_\xi 2^{-ns} \to 0$ as $n\to \infty$. 

(ii) By the $L^1_{\rm loc}$-convergence, we have
\[
  |u(x) - P\trace u(x)| \leq \bigg| u(x) - \dashint_{B(x)} T_{|x|} u \, d\mu\bigg| + \sum_{k \geq |x|}\bigg|\dashint_{B(x)}\big(T_{k+1}u - T_k u\big)\, d\mu\bigg|
\]
for all $x \in X$. The latter quantity can be estimated from above by a constant times
\[
  \sum_{\substack{|y| \geq |x|\\ B(y)\cap B(x) \neq \emptyset}} \frac{\mu\big(B(y)\big)}{\mu\big(B(x)\big)} |du(y)|;
\]
so Proposition \ref{pr:basic-operator} yields the desired conclusion.

(iii) This can easily be verified by using the density of Lipschitz functions in $L^1(Z)$. We refer to
\cite[Lemma 4.2]{BS} for details.
\end{proof}

Finally, we give the definition of our Triebel-Lizorkin space. Naturally, the definition could have been given much earlier---the auxiliary results proven so far are not needed for this purpose.

\defin{de:fspace}
Let $s \in (0,\infty)$, $p \in (Q/(Q+s),\infty)$, and $q \in (Q/(Q+s),\infty]$. Then the Triebel-Lizorkin space $\dfp$ is the vector space of all functions $f \in L^1_{\rm loc}(Z)$ such that
\[
 \|f\|_{\dfp} := \|d(Pf)\|_{\jp} < \infty.
\]
\edefin

\begin{remark} Recall that $d(Pf)$ is actually a vector-valued function on $X$ taking values in $\ell^2_{\degree}$. In the definition above, $\|d(Pf)\|_{\jp}$ stands for $\Vert \, |d(Pf)|\, \|_{\jp}$, where $|d(Pf)|=\|d(Pf)\|_{\ell^2_{\degree}}$ as before. We shall abuse notation in this way throughout the paper.
\end{remark}

The set $\dfp$ evidently becomes a quasi-normed space (with respect to the quasi-norm above) after dividing out the functions that are constant $\mu$-almost everywhere. We shall frequently abuse notation by writing $\dfp$ for both this quasi-normed space as well as the vector space of functions described above.

We have restricted the parameter $p$ in the above definition to be strictly greater than $Q/(Q+s)$, because our definition requires a priori local integrability and it is for these values of $p$ that we know the elements of the Fourier-analytically defined Triebel-Lizorkin spaces on an Euclidean space to be locally integrable.

\begin{remark} (i)\;
While the space $\jp$ depends on the exact choice of the hyperbolic filling $X$, the space $\dfp$ does not. This could be shown directly by a maximal function argument similar to the proof of Proposition \ref{pr:basic-operator}, but we are mostly interested in the case $0 < s \leq 1$, and in this case it is also a direct consequence of Proposition \ref{pr:identification-2} below. In fact, by examining the proof of Proposition \ref{pr:identification-2}, it can be seen that for given $s \in (0,1]$ and admissible values of $p$ and $q$, any two choices of $X$ yield equivalent quasi-norms on $\dfp$, with the equivalence constants independent of these two choices.

\smallskip
(ii) \; 
For some applications (see \cite{SS}), it is also useful to note that our methods allow some flexibility in the choice of the balls associated with the hyperbolic filling. More precisely, the parameters could be chosen so that $(\xi_x)_{x \in X_n}$ is for all $n$ a set of points in $Z$ with pairwise distances $\geq c_1 2^{-n}$ for some fixed constant $c_1$ (independent of $n$), that the radii $r_x$ corresponding to the balls $B(x) := B(\xi_x,r_x)$ ($x \in X_n$) are comparable to $2^{-n}$ uniformly in $x$ and $n$, and that the balls $\big(c_2 B(x)\big)_{x \in X_n}$ cover $Z$ for all $n$ where $c_2 \in (0,1)$ is a fixed constant. Under these assumptions, all the theorems of the present paper remain true.

\end{remark}

The results we have proven so far yield the following identification of $\dfp$ with a space of sequences on $X$. Below, we let $\jzero$ be the vector space of sequences $u$ in $\jp$ such that $u(\bzero) = 0$, and $\differencespace$ be the quasi-normed space of sequences $u$ on $X$ such that $|du| \in \jp$ and $u(\bzero) = 0$.

\thm{th:identification-1}
Let $0 < s < \infty$, $Q/(Q+s) < p < \infty$, and $Q/(Q+s) < q \leq \infty$.

\smallskip
\textnormal{(i)}
The trace $\trace u$ of a sequence $u$ in $\differencespace$ is zero $\mu$-almost everywhere if and only if $u \in \jzero$.

\smallskip
\textnormal{(ii)} We have an isomorphism 
\[
  \dfp \approx \differencespace / \jzero.
\]
of quasi-normed spaces. In particular, $\dfp$ is a quasi-Banach space, and when $1 < p,q < \infty$, it is a reflexive Banach space.
\ethm

\begin{proof}
(i) Assume first that $u \in \jzero$. Then $\trace u =0$ by Theorem \ref{th:trace-properties} (i). To prove the converse, let $u\in \differencespace$ with $\trace u=0$. In this case, $P\trace u =0$ so that Theorem \ref{th:trace-properties} (ii) implies $u\in \jp$. Since $u(\bzero)=0$, we have $u\in\jzero$.

(ii) We begin by verifying that $\differencespace$ is a quasi-Banach space in general, and a reflexive Banach space for $1 < p,\, q < \infty$. For the first statement, we argue just as in the proof of Proposition \ref{pr:quasi-banach} (i): $u_k(x) - u_k(x')$ is verified to converge for all $x$ and $x'$ joined by an edge in $E$ as $k \to \infty$, whence it follows that $u_k(x)$ converges for all $x$ as $k \to \infty$ because of the assumption that $u_k(\bzero) = 0$ for all $k$. In order to prove the reflexivity for $1 < p,q < \infty$, we note that in this case the space $\J^s_{p,q}(X, \ell^2_{\degree})$ of vector-valued sequences on $X$ (defined in the obvious manner) is a reflexive Banach space. Again this claim is verified exactly as the second statement in Proposition \ref{pr:quasi-banach} (i). Finally, $u \mapsto du$ is an isometry of $\differencespace$ onto a closed, and hence reflexive, subspace of $\J^s_{p,q}(X, \ell^2_\degree)$.

To proceed towards the claim concerning the isomorphism, we first note that $\jzero$ is a closed subspace of $\differencespace$. This follows by observing that the operator $\trace : \differencespace\to L^1(Z)$ is continuous (see Lemma \ref{le:trace}) and noting that by part (i) we have $\jzero=\trace^{-1}(\{ 0\})$. Hence $\differencespace / \jzero$ is well-defined. We claim that the map $H \: \dfp \to \differencespace / \jzero$, where
\[
  H(f):=Pf\; +\; \jzero\quad\text{for $f\in \dfp$}
\]
yields the desired isomorphism. Obviously, $H$ is well-defined and continuous by definition. The inverse map is given by
\[
  \trace_*(u+\jzero):=\trace u\quad \textrm{for $u\in \differencespace$.}
\]

In order to verify that $\trace_*\colon\differencespace / \jzero\to \dfp$ is continuous, we fix $$u+\jzero\in \differencespace / \jzero, $$ where $u$ is chosen so that
\[
  \| u +\jzero\|_{\differencespace / \jzero}\approx \| u\|_{\differencespace},
\]
and estimate 
\begin{align*}
    \| \trace_*(u +\jzero)\|_{\dfp} \;
& = \;\|\trace u\|_{\dfp}=\|d(P\trace u)\|_{\jp}\\
& \lesssim \;\|d(P\trace u-Pu)\|_{\jp}+\|d(Pu)\|_{\jp}\\
& \lesssim \;\|P\trace u-Pu\|_{\jp}+\|d(Pu)\|_{\jp}\\
& \lesssim \;\|d(Pu)\|_{\jp}\approx \|u+\jzero\|_{\differencespace / \jzero}.
\end{align*}
Here we applied Theorem \ref{th:trace-properties} (ii) together with the estimate $\| dg\|_{\jp}\lesssim \|g\|_{\jp}$, which follows from the fact that the graph $(X,E)$ has uniformly bounded valency.
\end{proof}

\section{Equivalence of definitions and density of Lipschitz functions}\label{se:equivalence}

We now turn to the identification of $\dfp$ with the Triebel-Lizorkin space $\dmp$ introduced by Koskela, Yang, and Zhou in \cite{KYZ}. We refer to \cite{GKZ} for a number of other characterizations of $\dmp$. The assumptions on $Z := (Z,d,\mu)$ here are the same as in the previous section.

Before stating the result, let us recall the definition of $\dmp$, where $0 < s < \infty$, $0 < p < \infty$ and $0 < q \leq \infty$. It is the vector space of measurable functions $f \: Z \to \complex$ such that
\[
  \|f\|_{\dmp} := \inf_{\vec{g} \in \D^s(f)} \bigg( \int_Z \big\|\big\{ g_k(\xi) : k \in \integer \big\}\big\|_{\ell^q}^p \,d\mu(\xi)\bigg)^{1/p}
\]
is finite, where the infimum is taken over the collection $\D^s(f)$ of all \emph{fractional $s$-Haj\l asz gradients} of $f$, i.e., sequences $\vec{g} := (g_k)_{k \in \integer}$ of measurable functions $g_k \: Z \to [0,\infty]$ such that there exists a measurable $E \subset Z$ with $\mu(E) = 0$ and
\[
  |f(\xi) - f(\eta)| \leq d(\xi,\eta)^s \big(g_k(\xi) + g_k(\eta)\big)
\]
whenever $\xi$, $\eta \in Z \backslash E$ and $2^{-k-1} \leq d(\xi,\eta) < 2^{-k}$. Note that this condition is void for $k$ such that $2^{-k-1}>\diam (Z)$. However, the above definition of $\dmp$ is also valid in the case when $Z$ is unbounded.

The set $\dmp$ becomes a quasi-normed space in an obvious way after dividing out the functions that are constant $\mu$-almost everywhere, and it is known that in case $Z = \real^d$, this quasi-normed space coincides with the standard Fourier-analytically defined Triebel-Lizorkin space $\dot{F}^s_{p,q}(\real^d)$ when $0 < s < 1$ and $p$, $q > d/(d+s)$. We refer to \cite{KYZ} for details.

\prop{pr:identification-2}
Let $0 < s \leq 1$, $Q/(Q+s) < p < \infty$, and $Q/(Q+s) < q \leq \infty$. Then
\[
  \dfp = \dmp
\]
with continuous embeddings.
\eprop

We remark here that since the smoothness index $s = 1$ is allowed, the space $\df^{1}_{p,\infty}(Z)$ coincides with the Haj\l asz space $\dot{M}^{1,p}(Z)$ for $Q/(Q+1) < p < \infty$ \cite[Proposition 2.1]{KYZ}. Recall that in the Euclidean setting, $\dot{M}^{1,p}(\real^d)$ coincides with the classical homogeneous Sobolev space for $1 < p < \infty$ \cite{H2}, and with a homogeneous Hardy-Sobolev space for $d/(d+1) < p \leq 1$ \cite{KS}. On the other hand, Proposition \ref{pr:identification-2} above together with \cite[Theorem 4.1]{GKZ} implies that if $Z$ supports a weak $(1,p)$-Poincar\'e inequality for some $p > 1$, then $\df^1_{p,q}(Z)$ is trivial (i.e., only contains constant functions) for all $q \in (Q/(Q+1),\infty)$.

For the proof of Proposition \ref{pr:identification-2}, we need the following Poincar\'e-type inequality. In the Euclidean setting it is a special case of \cite[Lemma 2.3]{KYZ}, and for doubling metric measure spaces as in our setting, it can be proven with a similar argument as noted in the proof of \cite[Theorem 4.1]{KYZ}. It readily implies that the elements of $\dmp$ belong to $L^1_{\rm loc}(Z)$ when $0 < s \leq 1$ and $p > Q/(Q+s)$.

\lem{le:poincare}
Let $s \in (0,1]$ and $0 < \epsilon < \epsilon' < s$. Then for all $\xi \in Z$, $k \in \integer$, measurable functions $f$ on $Z$, and $\vec{g} = (g_j)_{j\in\integer} \in {\mathbb D}^s(f)$ we have
\[
  \inf_{c\in\complex}\, \dashint_{B(\xi,2^{-k})} |f - c| \,d\mu \leq C 2^{-k\epsilon'} \sum_{j\geq k-2} 2^{-j(s-\epsilon')} \bigg( \dashint_{B(\xi,2^{-k+1})}g_j^{\frac{Q}{Q+\epsilon}} \,d\mu\bigg)^{\frac{Q+\epsilon}{Q}},
\]
with a constant $C$ independent of $\xi$, $k$, $f$, and $\vec g$.
\elem

\begin{proof}[Proof of Proposition \ref{pr:identification-2}]
To simplify the notation, we only consider the case $q < \infty$, as the case $q = \infty$ is easier. To establish the embedding $\dmp \subset \dfp$, let $f \in \dmp$ and $\vec g := (g_j)_{j \in \integer} \in {\mathbb D}^s(f)$. Then $f \in L^1_{\rm loc}(Z)$; so if $x$ and $x'$ are distinct points of $X$ joined by an edge in $E$, we can find a ball $B$ that has a radius comparable to $2^{-|x|}$ and covers both $B(x)$ and $B(x')$.
Setting 
\[
  f_B:=\dashint_{B} f \,d\mu,
\]
we obtain 
\begin{align*}
  |Pf(x) - Pf(x')| &\leq \dashint_{B(x)}|f - f_B| \,d\mu + \dashint_{B(x')} |f - f_B| \, d\mu \\& \lesssim \dashint_{B}|f -f_B| \,d\mu.
\end{align*}
If we take $\vec g := (g_j)_{j\in\integer} \in {\mathbb D}^s(u)$ and $\epsilon$, $\epsilon' \in (0,s)$ so that $\epsilon < \epsilon'$ and $p$, $q > Q/(Q+\epsilon)$, Lemma \ref{le:poincare} gives
\[
  |Pf(x)-Pf(x')| \lesssim 2^{-|x|\epsilon'} \sum_{j \geq |x|-\sigma} 2^{-j(s-\epsilon')} \Bigg( \dashint_{\sigma B(x)}g_j^{\frac{Q}{Q+\epsilon}} \,d\mu\Bigg)^{\frac{Q+\epsilon}{Q}}
\]
for some uniform constant $\sigma \geq 1$. Since the graph $(X,E)$ has bounded valency, the left-hand side of the estimate above can be replaced by $|d(Pf)(x)|$. In particular, if $x \in \gox$ for some $\xi \in Z$, we have
\[
  |d(Pf)(x)| \lesssim 2^{-|x|\epsilon'} \sum_{j \geq |x|-\sigma} 2^{-j(s-\epsilon')} \hlmax \Big( g_j^{\frac{Q}{Q+\epsilon}} \Big)(\xi)^{\frac{Q+\epsilon}{Q}}.
\]
Using the subadditivity of $t \mapsto t^q$ or H\"older's inequality depending on $q$, we get
\[
  \sum_{x\in\gox} \big[2^{|x|s} |d(Pf)(x)|\big]^q \lesssim \sum_{x\in\gox} 2^{(s-\epsilon')(q\land 1)|x|} \sum_{j \geq |x|-\sigma} 2^{-j(s-\epsilon')(q\land 1)} \hlmax \Big( g_j^{\frac{Q}{Q+\epsilon}} \Big)(\xi)^{\frac{Q+\epsilon}{Q}q},
\]
and by changing the order of summation and using again the fact that $\# (\gox \cap X_j) \approx 1$, we have
\[
  \sum_{x\in\gox} \big[2^{|x|s} |d(Pf)(x)|\big]^q \lesssim \sum_{j\in\integer} \hlmax \Big( g_j^{\frac{Q}{Q+\epsilon}} \Big)(\xi)^{\frac{Q+\epsilon}{Q}q}.
\]
Since $\frac{Q+\epsilon}{Q}p$, $\frac{Q+\epsilon}{Q}q > 1$, the Fefferman-Stein maximal theorem then gives
\[
  \|f\|_{\dfp}^{p} \lesssim \int_Z \bigg( \sum_{j\in\integer} \hlmax \Big( g_j^{\frac{Q}{Q+\epsilon}} \Big)(\xi)^{\frac{Q+\epsilon}{Q}q}\bigg)^{ \frac{\frac{Q+\epsilon}{Q}p}{\frac{Q+\epsilon}{Q}q} } \,d\mu(\xi) \lesssim \int_Z \Big(\sum_{j\in\integer} g_j(\xi)^q\Big)^{p/q} \,d\mu(\xi).
\]
If we take the infimum over $\vec g \in {\mathbb D}^s(f)$, we finally obtain $\|f\|_{\dfp} \lesssim \|f\|_{\dmp}$.

To establish the embedding $\dfp \subset \dmp$, let $f \in \dfp$ and $\xi_1$ and $\xi_2$ be distinct Lebesgue points of $f$. Denote by $k$ the unique integer with $2^{-k-1} \leq d(\xi_1,\xi_2) < 2^{-k}$. Then $k\geq -1$. For $i \in \{1,2\}$ and $n \geq (k-1)_+:=\max(k-1,0)$ we can then choose points $x^i_n \in X_n$ so that $\xi_i \in B(x^i_n)$. Since the balls $\frac12 B(x)$, $x \in X_{(k-1)_+}$, cover $Z$, we can further assume that $x^1_{(k-1)_+} = x^2_{(k-1)_+}$. We thus get

\begin{align*}
  |f(\xi_1) - f(\xi_2)|\; & = \;\bigg| \sum_{n \geq (k-1)_+} \Big( Pf(x^1_{n+1}) - Pf(x^1_{n}) \Big) - \sum_{n \geq (k-1)_+} \Big( Pf(x^2_{n+1}) - Pf(x^2_{n}) \Big) \bigg|\\
  & \leq \;\sum_{\substack{x \in \Gamma_{1,\xi_1} \\ |x| \geq k}} |d(Pf)(x)| + \sum_{\substack{x \in \Gamma_{1,\xi_2} \\ |x| \geq k}} |d(Pf)(x)|,
\end{align*}
so that
\[
  |f(\xi_1) - f(\xi_2)| \lesssim d(\xi_1,\xi_2)^s\Bigg( 2^{ks} \sum_{\substack{x \in \Gamma_{1,\xi_1} \\ |x| \geq k}} |d(Pf)(x)| + 2^{ks}\sum_{\substack{x \in \Gamma_{1,\xi_2} \\ |x| \geq k}} |d(Pf)(x)| \Bigg).
\]

Since the set of Lebesgue points of $f$ has full $\mu$-measure in $Z$, the sequence of functions $(g_k)_{k\in\integer}$ defined by
\[
  g_k(\xi) := 2^{ks} \sum_{\substack{x \in \gox \\ |x| \geq k}} |d(Pf)(x)|
\]
for all $k \in \integer$ is therefore a constant times an element of ${\mathbb D}^s(f)$. As before we have
\[
  \sum_{k\in\integer} g_k(\xi)^q \lesssim \sum_{k \in \integer} 2^{ks(q\land 1)} \sum_{\substack{x \in \gox \\ |x| \geq k}} 2^{-|x|s(q\land 1)} \big[2^{|x|s}|d(Pf)(x)|\big]^q \approx \sum_{x\in\gox} \big[2^{|x|s}|d(Pf)(x)|\big]^q, 
\]
so that $\|f\|_{\dmp} \lesssim \|f\|_{\dfp}$.
\end{proof}

The trace approximations $T_n(Pf)$ defined in Section \ref{se:definitions} provide a natural discrete convolution type approximation for our functions $f\in\dfp$, as Theorem \ref{th:trace-properties} shows that they converge to $f$ in $L^1(Z)$. We next verify that when $q < \infty$, then we actually have convergence with respect to the Triebel-Lizorkin norm.

\thm{th:approximations}
\textnormal{(i)} Let $0 < s < 1$, $Q/(Q+s) < p < \infty$, and $Q/(Q+s) < q \leq \infty$, or $s = 1$, $Q/(Q+1) < p < \infty$, and $q = \infty$. Then
\[
  \|T_n (Pf)\|_{\dfp} \lesssim \|f\|_{\dfp}
\]
for all $f \in \dfp$ and $n \in \nanu$, with an implicit constant independent of $f$ and $n$.

\smallskip
\textnormal{(ii)} Let $0 < s < 1$, $Q/(Q+s) < p,q < \infty$, and $f \in \dfp$. Then $T_n(Pf) \to f$ in $\dfp$ as $n \to \infty$.
\ethm

\begin{proof} (i) Fix $n\geq 1$. We first estimate
  $|dP(T_nPf)(x)|$
in case $|x| \leq n+2$. Let $x'$ be an element of $X$ joined with $x$ by an edge in $E$. We have 
\begin{align*}
  |P(T_nPf)(x) - P(T_nPf)(x')| \;
& \leq \;|d(Pf)(x)| + | P(T_nPf)(x) - Pf(x) | \\
& \quad \quad + | P(T_nPf)(x') - Pf(x') |.
\end{align*}

As $T_k(Pf) \to f$ pointwise $\mu$-almost everywhere as $k\to\infty$, we have (compare with the proof of Lemma \ref{le:trace})
\begin{align*}
  | P(T_n(Pf))(x) - Pf(x) | \;
& \leq \; \dashint_{B(x)} \big|f - T_n(Pf)\big| \,d\mu \\
& \leq \;\sum_{k \geq n} \dashint_{B(x)} \big| T_{k+1}(Pf) - T_k(Pf) \big| d\mu \\
& \lesssim \;\sum_{\substack{|y| \geq n \\ B(y) \cap B(x) \neq \emptyset}} \frac{\mu\big(B(y)\big)}{\mu\big(B(x)\big)}|d(Pf)(y)| \\
& = \;\sum_{\substack{|y| \geq |x|-2 \\ B(y) \cap B(x) \neq \emptyset}} \frac{\mu\big(B(y)\big)}{\mu\big(B(x)\big)}|d(Pf)(y)| \chi_{\cup_{k \geq n}X_k}(y).
\end{align*}

In particular, taking $x_* \in X$ so that $|x_*| = (|x|-2)_+$ and $B(x) \subset B(x_*)$, the last  quantity can be estimated from above by a constant times $T(|d(Pf)|\chi_{\cup_{k\geq n}X_k})(x_*)$, with $T$ as in Proposition \ref{pr:basic-operator}. The term $| P(T_n(Pf))(x') - Pf(x') |$ above can be estimated in a similar manner. We thus have
\beqla{eq:approximations-1}
  |dP(T_nPf)(x)| \lesssim |d(Pf)(x)| + T(|d(Pf)|\chi_{\cup_{k \geq n}X_k})(x_*) \quad \text{when $|x| \leq n+2$}.
\eeq

On the other hand, if $x \in \gox$ for some $\xi \in Z$ and $|x| \geq n+3$, we may for all neighbors $x'$ of $x$ take $y_{x,x'} \in X_n$ so that $B(y_{x,x'})$ covers both $B(x)$ and $B(x')$. We estimate $|P(T_n(Pf))(x) - P(T_n(Pf))(x')|$ from above by
\[
  \dashint\dashint_{B(x) \times B(x')} \sum_{\substack{y \in X_n \\ B(y) \cap (B(x) \cup B(x')) \neq \emptyset}}|Pf(y)-Pf(y_{x,x'})||\psi_y(\eta) - \psi_y(\eta')|\, d\mu(\eta) \, d\mu(\eta').
\]

Using the Lipschitz continuity of the functions $\psi_y$, we get
\[
  |dP(T_n(Pf))(x)| \lesssim 2^{n-|x|} \sum_{\substack{y \in \gox \\ |y| = n}} |d(Pf)(y)|,
\]
so that
\beqla{eq:approximations-2}
\sum_{\substack{x \in \gox \\ |x| \geq n+3}} \big[2^{|x|s}|dP(T_nPf)(x)|\big]^q \lesssim \sum_{\substack{y \in \gox \\ |y| = n}} \big[2^{|y|s}|d(Pf)(y)|\big]^q
\eeq
(with an obvious modification for $q = \infty$). Combining \eqref{eq:approximations-1} and \eqref{eq:approximations-2} with Proposition \ref{pr:basic-operator}, we obtain the desired upper bound for $\|T_n(Pf)\|_{\dfp}$.

(ii) This is almost contained in the argument above. We have
\beqla{eq:approximations-3}
  |dP(f - T_n(Pf))(x)| \lesssim T(|d(Pf)| \chi_{\cup_{k \geq n}X_k})(x_*) \quad \text{when} \quad |x| \leq n+2,
\eeq
with $x_*$ is as before. If, on the other hand, $x \in \gox$ for some $\xi$ and $|x| \geq n+3$, we have
\begin{align}
  |dP(f - T_nPf)(x)| \;
& \lesssim \;|d(Pf)(x)| + |d(P(T_n(Pf)))(x)| \notag\\
& \lesssim \;|d(Pf)(x)| + 2^{n-|x|} \sum_{\substack{y \in \gox \\ |y| = n}} |d(Pf)(y)|\label{eq:approximations-4}. 
\end{align}

Combining \eqref{eq:approximations-3} and \eqref{eq:approximations-4} with Proposition \ref{pr:basic-operator}, we obtain
\[
  \|f - T_n(Pf)\|_{\dfp} \lesssim \Big( \int_{Z} \big\| \{ 2^{|x|s} |d(Pf)(x)| : x \in \gox,\, |x| \geq n \} \big\|_{\ell^q}^p \,d\mu(\xi)\Big)^{1/p},
\]
where the latter quantity tends to zero as $n \to \infty$ by the dominated convergence theorem.
\end{proof}

The theorem above immediately yields the corollary below in case that the space $Z$ is bounded. If $Z$ is unbounded, we need some additional observations, which follow from Proposition \ref{pr:retraction}. We postpone the details until the end of Section \ref{se:interpolation}. 

\cor{co:lipschitz-density} Suppose that $0 < s < 1$ and $Q/(Q+s) < p,\,q < \infty$. Then Lipschitz functions with bounded support are dense in $\dfp$. 
\ecor

In the setting of an unbounded metric measure space that supports a ``reverse doubling'' property in addition to the doubling property \eqref{eq:doubling}, a result similar to the Corollary above can be found in \cite[Proposition 5.21]{HMY}.

\section{Unbounded spaces and embeddings}\label{se:unbounded}

We have so far assumed the metric space $Z$ to be bounded, but this is actually not an essential assumption. For an unbounded space $Z$, we simply choose a maximal set of points $\{\xi_x : x \in X_n\} \subset Z$ with pairwise distances at least $2^{-n-1}$ for all $n \in \integer$ (not just positive $n$), set $X := \sqcup_{n \in \integer } X_n$ and write $|x| = n$ for all $x \in X_n$. The Poisson extension $Pf$ is now well-defined for any $f \in L^1_{\rm loc}(Z)$. Definition \ref{de:sspace} and Proposition \ref{pr:quasi-banach} generalize without any changes. Only in Lemma \ref{le:trace} (and consequently in the results that depend on this lemma) do we use the fact that $Z$ is bounded; in general, Triebel-Lizorkin functions of course are not integrable over an unbounded space.

In the unbounded case, the limit $\trace u := \lim_{n \to \infty} T_n u$ exists in $L^1_{\rm loc}(Z)$ and pointwise almost everywhere whenever $|du| \in \jp$, and the embedding in Lemma \ref{le:trace} may be replaced with
\[
  \|\trace u - u(x_B)\|_{L^1(B)} \lesssim c(r,\mu(B),p,s)\|du\|_{\jp},
\]
where $B$ is any ball of radius $r$ and $x_B$ is any such element of $X$ that $4r \leq 2^{-|x_B|} < 8r$ and $B(x_B)$ contains $B$ as a subset. The implicit constant in the previous inequality is independent of $u$, $B$, and $x_B$. In particular, if we fix a ball $B_0$ of radius $1$ and $x_0 \in X_{-2}$ corresponding to $B_0$ as in the embedding above, we may define $\differencespace$ as the space of sequences $u$ on $X$ such that $|du| \in \jp$ and $u(x_0) = 0$. The trace operator $\trace$ then takes $\differencespace$ continuously into $L^1_{\rm loc}(Z)$ endowed with the topology induced by the seminorms $f \mapsto \|f\|_{L^1(B_0)}$ and $f \mapsto \inf_{c\in\complex}\|f-c\|_{L^1(kB_0)}$, $k \geq 2$. The analogs to the rest of our results can then be proven with slight modifications. In the rest of the paper, we state our results for general $Z$ that can be bounded or unbounded.

In the remainder of this section we establish general embedding results for the Triebel-Lizorkin spaces. Suppose momentarily that $Z$ is unbounded and that the measure $\mu$ is $Q$-Ahlfors regular, i.e., $\mu(B(\xi,r)) \approx r^Q$ uniformly in $\xi \in Z$. Fix $\xi_0 \in Z$ and let $f_\lambda$ for $\lambda > 0$ stand for the function $\xi \mapsto (1-\lambda d(\xi,\xi_0))_+$. A straightforward computation then yields
\beqla{eq:bump-functions}
  \|f_\lambda\|_{\dfp} \approx \lambda^{s - Q/p}
\eeq
when $0 < s < 1$, $Q/(Q+s) < p < \infty$ and $Q/(Q+s) < q \leq \infty$, with the implicit constants independent of $\lambda$. In particular, to have an embedding of the type $\df^{s_0}_{p_0,q_0}(Z) \subset \df^{s_1}_{p_1,q_1}(Z)$, we need to have $s_0 - Q/p_0 = s_1 - Q/p_1$.

The proposition below establishes all such embeddings with $p_1 > p_0$ under the slightly weaker assumption that $\mu$ supports the lower mass bound $\mu(B(\xi,r)) \geq c_0 r^Q$ for all $\xi \in Z$ and $0 < r < \diam Z$. If $Z$ is bounded, this follows automatically from the doubling condition on $\mu$. In the Euclidean case, these embeddings are a special case of a result due to Jawerth \cite[Theorem 2.1 (ii)]{J}, and the proof below is based on his methods. We apply the result in the proof of Proposition \ref{pr:qs-converse} below.

\prop{pr:jawerth-embedding}
Suppose that $\mu$ satisfies the doubling condition \eqref{eq:doubling} and supports the lower mass bound $\mu(B(\xi,r)) \geq c_0r^Q > 0$ for all $\xi \in Z$ and $0 < r < \diam Z$. Let $0 < s_1 < s_0 < \infty$, $0 < p_0 < p_1 < \infty$ and $0 < q_0,q_1 \leq \infty$ be such that $s_0 - Q/p_0 = s_1 - Q/p_1$. Then
\[
  \J^{s_0}_{p_0,q_0}(X) \subset \J^{s_1}_{p_1,q_1}(X)
\]
with a continuous embedding. In particular, if $Q/(Q+s_i) < \min(p_i,q_i)$ in addition to the above assumptions, then
\beqla{eq:jawerth-embedding}
  \df^{s_0}_{p_0,q_0}(Z) \subset \df^{s_1}_{p_1,q_1}(Z)
\eeq
with a continuous embedding.
\eprop

\begin{proof}
It suffices to consider the case with $q_0 = \infty$ and $q_1 < \infty$. Let $u \in \J^{s_0}_{p_0,\infty}(X)$. Let $t > 0$ be given and suppose that $\xi \in Z$ is a point such that
\[
  \bigg(\sum_{x \in \gox} \big[2^{|x|s_1}|u(x)|\big]^{q_1}\bigg)^{1/q_1} > t.
\]
Now, since 
\[
  |u(x)| \leq \mu(B(x))^{-1/p_0} 2^{-|x|s_0}\|u\|_{\J^{s_0,(1)}_{p_0,\infty}(X)} \leq c_1 2^{|x|(Q/p_1 - s_1)}\|u\|_{\J^{s_0}_{p_0,\infty}(X)}
\]
for all $x \in X$ with a constant $c_1 > 0$, we have
\[
  \bigg(\sum_{\substack{x \in \gox\\ |x| \leq N}} \big[2^{|x|s_1}|u(x)|\big]^{q_1}\bigg)^{1/q_1} \leq c_2 2^{Q N/p_1} \|u\|_{\J^{s_0}_{p_0,\infty}(X)}
\]
for all integers $N$ with another constant $c_2 > 0$. In particular, we take $N$ so that $c_3 t \leq 2^{QN/p_1}\|u\|_{\J^{s_0}_{p_0,\infty}(X)} < c_4t$, where the constants $c_3$ and $c_4$ are chosen so that the estimate above guarantees
\[
  \bigg(\sum_{\substack{x \in \gox\\ |x| > N}} \big[2^{|x|s_1}|u(x)|\big]^{q_1}\bigg)^{1/q_1} \gtrsim t.
\]
Now since $s_1 < s_0$, we have
\begin{align*}
  \bigg(\sum_{\substack{x \in \gox\\ |x| > N}} \big[2^{|x|s_1}|u(x)|\big]^{q_1}\bigg)^{1/q_1} & \lesssim 2^{N(s_1-s_0)} \sup_{x \in \gox} 2^{|x|s_0}|u(x)| = 2^{N(Q/p_1 - Q/p_0)} \sup_{x \in \gox} 2^{|x|s_0}|u(x)| \\
 & \approx t^{1 - p_1/p_0} \|u\|_{\J^{s_0}_{p_0,\infty}(X)}^{p_1/p_0 - 1} \sup_{x \in \gox} 2^{|x|s_0}|u(x)|,
\end{align*}
so that
\[
  \sup_{x \in \gox} 2^{|x|s_0}|u(x)| > c_5 t^{p_1/p_0} \|u\|_{\J^{s_0}_{p_0,\infty}(X)}^{1 - p_1/p_0}
\]
for some constant $c_5 > 0$. In particular,
\begin{align*}
  \|u\|_{\J^{s_1}_{p_1,q_1}(X)}^{p_1} \;
& \approx \;\int_{0}^{\infty} t^{p_1 - 1} \mu\bigg(\bigg\{\xi : \bigg(\sum_{x \in \gox} \big[2^{|x|s_1}|u(x)|\big]^{q_1}\bigg)^{1/q_1} > t \bigg\} \bigg) \,dt \\
& \leq \;\int_{0}^{\infty} t^{p_1 - 1} \mu\bigg(\bigg\{\xi : \sup_{x \in \gox} 2^{|x|s_0}|u(x)| > c_5 t^{p_1/p_0} \|u\|_{\J^{s_0}_{p_0,\infty}(X)}^{1 - p_1/p_0} \bigg\} \bigg) dt \\
& \approx \;\|u\|_{\J^{s_0}_{p_0,\infty}(X)}^{p_1-p_0} \int_{0}^{\infty} t^{p_0 - 1} \mu\bigg(\bigg\{\xi : \sup_{x \in \gox} 2^{|x|s_0}|u(x)| > t \bigg\} \bigg) \,dt \\
& \approx \;\|u\|_{\J^{s_0}_{p_0,\infty}(X)}^{p_1}. \qedhere
\end{align*}
\end{proof}

As $\dfp \subset \df^{s}_{p,\infty}(Z) = \dm^s_{p,\infty}(Z)$ with continuous embeddings for $0 < s \leq 1$ and admissible values of $p$ and $q$, we may essentially replace the right-hand side of \eqref{eq:jawerth-embedding} with $L^{p_{*}}(Z)$, where $p_* = Qp_0/(Q-p_0s_0)$, for $0 < s_0 \leq 1$ and $Q/(Q+s_0) < p_0 < Q/s_0$. We refer to \cite[Lemma 4.2]{KYZ2} and the references therein for details.

Back in our original setting, where $Z$ is only assumed to satisfy the doubling condition \eqref{eq:doubling}, we do not have global estimates for the measures of balls, but the argument in the proof of Proposition \ref{pr:jawerth-embedding} may be localized in order to obtain
\begin{align*}
  & r^{s_1}\bigg( \dashint_{B} \big\| \{ 2^{|x|s_1} |u(x)| : x \in \gox,\; 2^{-|x|} < r/2 \} \big\|_{\ell^{q_1}}^{p_1} \,d\mu(\xi)\bigg)^{1/p_1} \\
  & \qquad \lesssim r^{s_0}\bigg( \dashint_{2B} \big\| \{ 2^{|x|s_0} |u(x)| : x \in \gox,\; 2^{-|x|} < r/2 \} \big\|_{\ell^{q_0}}^{p_0} \,d\mu(\xi)\bigg)^{1/p_0}
\end{align*}
for all balls $B$ with radius $r > 0$, with the implicit constant independent of $u$ and $B$. We omit the details.

\section{Quasisymmetric invariance}\label{se:invariance}

Here we present a simple proof of one of the main results of \cite{KYZ}, namely that the composition operator $f \mapsto f\circ \varphi$ induced by a quasisymmetric homeomorphism $\varphi \: Z \to Z'$, where $Z := (Z,d,\mu)$ and $Z' := (Z,d',\mu')$ are $Q$-Ahlfors regular and have sufficiently reasonable geometry, is bounded from $\df^{s}_{Q/s,q}(Z')$ to $\df^{s}_{Q/s,q}(Z)$ whenever $0 < s \leq 1$ and $Q/(Q+s) < q \leq \infty$. Note that, by \eqref{eq:bump-functions}, we can in general expect the function space $\df^{s}_{p,q}$ to have some sort of conformal invariance only when $p = Q/s$. All measures appearing in this chapter are assumed to be Borel regular.

Let us first recall the definition of quasisymmetry and formulate the assumptions on the metric measure spaces under consideration. The quasisymmetry assumption on the homeomorphism $\varphi \: Z \to Z'$ means that there exists an increasing homeomorphism $\rho \: (0,\infty) \to (0,\infty)$ so that
\[
  \frac{d'\big(\varphi(\xi),\varphi(\eta_1)\big)}{d'\big(\varphi(\xi),\varphi(\eta_2)\big)} \leq \rho\bigg(\frac{d(\xi,\eta_1)}{d(\xi,\eta_2)} \bigg)
\]
for all $\xi \in Z$ and $\eta_1$, $\eta_2 \in Z \backslash \{\xi\}$. It follows easily that then $\varphi^{-1} \: Z' \to Z$ is also quasisymmetric. In Proposition \ref{pr:j-qsinvariance} below, we assume that $Z$ and $Z'$ are both locally compact and $Q$-Ahlfors regular for some $Q > 1$, and that $Z$ is complete and supports a weak $(1,Q)$-Poincar\'e inequality. Note that the assumption concerning the Poincar\'e inequality implies that $Z$ is connected. We refer to e.g.~\cite{HK} for the necessary definitions, as well as examples of spaces with these properties. 

Our assumptions imply that the pullback measure $\sigma_\varphi := E \mapsto \mu'(\varphi(E))$ is doubling on $Z$ and that the ``Jacobian''
\[
  J_\varphi(\xi) := \lim_{r \to 0} \frac{ \mu'\big( \varphi\big( B_d(\xi,r) \big) \big)}{\mu\big(B_d(\xi,r)\big) }
\]
is well-defined for $\mu$-almost all $\xi \in Z$ and satisfies the standard change of variables formula for any integrable function (i.e.~$d\sigma_\varphi = J_\varphi \,d\mu$), as well as the reverse H\"older inequality
\[
  \bigg(\dashint_{B} J_\varphi^{R_\varphi} \,d\mu \bigg)^{1/R_\varphi} \leq C_\varphi \dashint_{B} J_\varphi \,d\mu
\]
with some exponent $R_\varphi > 1$ and a constant $C_\varphi \geq 1$ for all balls $B$ on $Z$. We refer to \cite{KZ}, \cite[Remark 3.4]{HK} as well as the references therein and \cite[Theorem 7.11]{HK} for details.

Our proof of the quasi-invariance of the Triebel-Lizorkin spaces (Theorem \ref{th:f-qsinvariance} below) is based on moving the action of $\varphi$ inside $X$. This is possible since it is actually known that $\varphi$ extends to a quasi-isometric mapping $\Phi$ between the hyperbolic fillings $X$ and $X'$. We next discuss the basic properties of this extension. First of all, $(X,E)$ can be viewed as a metric graph equipped with the natural length metric (so that every edge has length $1$). We denote by $d_X$ the metric on the set of vertices $X$ induced by this length metric, and define $d_{X'}$ on $X'$ in a similar manner. Then there exists a mapping $\Phi \: X \to X'$ that extends $\varphi$ in the sense that $\varphi(\xi) \in B(\Phi(x))$ whenever $\xi \in B(x)$. More explicitly, one may define $\Phi(x)$ as an element at a maximal level of $X'$ so that $\varphi(B(x)) \subset B(\Phi(x))$. The quasi-isometry property of $\Phi$ means that 

\beqla{eq:quasiisometry}
  \frac1\lambda d_X (x,y) - c \leq d_{X'}\big(\Phi(x),\Phi(y)\big) \leq \lambda d_X (x,y) + c
\eeq
for some constants $\lambda \geq 1$, $c \geq 0$ whenever $x$, $y \in X$ and
\[
  \bigcup_{x \in X} B_{d_{X'}}\big(\Phi(x),r\big) = X'
\]
for some $r > 0$. We refer to \cite[Theorem 7.2.1]{BuS} for details. 

The following two properties of $\Phi$ will be used in the proof below:
\beqla{eq:qi-property1}
  B\big( \Phi(x) \big) \subset \varphi\big(\sigma B(x)\big) \text{ for all } x \in X
\eeq
with $\sigma > 1$ independent of $x$ and
\beqla{eq:qi-property2}
  \#\big\{ x \in \gox : |\Phi(x)| = k \big\} \lesssim 1 \text{ for all } \xi \in Z,\; k \in \integer.
\eeq
The inclusion \eqref{eq:qi-property1} is a consequence of the definition of $\Phi$ and the quasisymmetry of $\varphi$, while \eqref{eq:qi-property2} follows easily from \eqref{eq:quasiisometry} and the fact that $$d_{X'}(\Phi(x),\Phi(y)) \leq ||\Phi(x)| - |\Phi(y)|| + 1$$ for all $x$, $y \in \gox$.

\prop{pr:j-qsinvariance}
Let the spaces $Z$ and $Z'$ be as in the discussion above, and let $0 < s < \infty$ and $0 < q \leq \infty$. Then the composition operator induced by $\Phi$ is bounded from $\J^{s}_{Q/s,q}(X')$ into $\J^{s}_{Q/s,q}(X)$, in the sense that
\[
  \| u\circ \Phi \|_{\J^{s}_{Q/s,q}(X)} \lesssim \| u \|_{\J^{s}_{Q/s,q}(X')}
\]
for all $u \in \J^{s}_{Q/s,q}(X')$, with the implicit constant independent of $u$.
\eprop

\begin{proof}
We begin by defining two auxiliary parameters as in the proof of \cite[Theorem 1.3]{KYZ}: take $\delta > 0$ and $p > 1$ so that
\beqla{eq:qsinvariance-params}
  \delta p < \min(q , Q/s) \quad \text{and} \quad 0 <\frac{Q - \delta s}{Q} \cdot \frac{p}{p - 1} < R_\varphi.
\eeq
One can simply choose $p$ large enough so that $p/(p-1)<R_\varphi$ and then take $\delta<Q/s$ small enough so that the second inequality in \eqref{eq:qsinvariance-params} holds.

Let $u \in \J^{s}_{Q/s,q}(X')$ and define the functions $U_k$ on $Z'$ for $k \in \integer$ by setting
\[
  U_k(\xi') := \| 2^{|x'|s} |u(x')| : x' \in \Gamma_{1,\xi'} \cap X'_k \|_{\ell^q}.
\]
Then
\[
  \| u \|_{\J^{s}_{Q/s,q}(X')}^{Q/s} = \int_{Z'} \bigg(\sum_{k \in \integer} U_k(\xi')^q \bigg)^{(Q/s)/q} \,d\mu'(\xi')
\]
(with an obvious modification for $q = \infty$). Now if $\xi \in Z$ and $x \in \gox$, we have
\beqla{eq:qsinvariance-estimate}
  2^{|x|s} |(u\circ \Phi)(x)| \leq 2^{(|x| - |\Phi(x)|)s} \bigg(\dashint_{B(\Phi(x))}U_{|\Phi(x)|}(\eta')^{\delta} \,d\mu'(\eta')\bigg)^{1/\delta},
\eeq
and by \eqref{eq:qi-property1}, there is $\sigma > 1$ so that the quantity inside the parentheses can be estimated from above by a constant times

\begin{align*}
& \qquad \qquad \frac{\mu\big( \sigma B(x))}{\mu'\big(B(\Phi(x))\big)} \dashint_{\sigma B(x)} U_{|\Phi(x)|}\big(\varphi(\eta)\big)^{\delta} J_\varphi(\eta)\, d\mu(\eta) \\
& \lesssim \;\frac{\mu\big(\sigma B(x))}{\mu'\big(B(\Phi(x))\big)} \bigg( \dashint_{\sigma B(x)} U_{|\Phi(x)|}\big(\varphi(\eta)\big)^{\delta p} J_\varphi(\eta)^{\delta p s/Q} \,d\mu(\eta)\bigg)^{1/p} \\
& \qquad \qquad \times \bigg( \dashint_{\sigma B(x)} J_\varphi(\eta)^{[(Q-\delta s)p]/[Q(p-1)]} \,d\mu(\eta)\bigg)^{(p-1)/p} \\
& \lesssim \;\frac{\mu\big(\sigma B(x))}{\mu'\big(B(\Phi(x))\big)} \hlmax \Big( [ U_{|\Phi(x)|} \circ \varphi \big]^{\delta p} J_\varphi^{\delta p s/Q}\Big)(\xi)^{1/p} \\
& \qquad \qquad \times \bigg[ \frac{\mu'\big(\varphi(\sigma B(x)))}{\mu\big(\sigma B(x)\big)} \bigg]^{1 - (\delta s )/Q} \\
& \lesssim \; 2^{-(|x| - |\Phi(x)|)\delta s} \hlmax \Big( [ U_{|\Phi(x)|} \circ \varphi \big]^{\delta p} J_\varphi^{\delta p s/Q}\Big)(\xi)^{1/p}.
\end{align*}
Here we used H\"older's inequality in the first inequality, the reverse H\"older inequality and the change of variables formula in the second inequality, and the Ahlfors regularity of both metric measure spaces in the third inequality. Substituting this into \eqref{eq:qsinvariance-estimate}, we obtain
\[
  2^{|x|s} |(u\circ \Phi)(x)| \lesssim \hlmax \Big( [ U_{|\Phi(x)|} \circ \varphi \big]^{\delta p} J_\varphi^{\delta p s/Q}\Big)(\xi)^{1/(\delta p)},
\]
and by \eqref{eq:qi-property2}, we get
\[
  \| 2^{|x|s} |(u\circ \Phi)(x)| : x \in \gox \|_{\ell^q} \lesssim \bigg(\sum_{k \in \integer} \hlmax \Big( [ U_{k} \circ \varphi \big]^{\delta p} J_\varphi^{\delta p s/Q}\Big)(\xi)^{q/(\delta p)} \bigg)^{1/q}
\]
(with an obvious modification for $q = \infty$). As $\delta p < \min(q , Q/s)$, the Fefferman-Stein maximal theorem and the change of variables formula thus yield
\begin{align*}
  \| u\circ \Phi \|_{\J^{s}_{Q/s,q}(X)}^{Q/s} \;
& \lesssim \;\int_{Z} \bigg( \sum_{k \in \integer} U_k\big( \varphi(\xi)\big)^q \bigg)^{(Q/s)/q} J_\varphi (\xi) \,d\mu(\xi) \\
& = \;\int_{Z'} \bigg( \sum_{k \in \integer} U_k(\xi')^q \bigg)^{(Q/s)/q} \,d\mu'(\xi') = \| u \|_{\J^{s}_{Q/s,q}(X')}^{Q/s}. \qedhere
\end{align*}
\end{proof}

The corresponding result for our Triebel-Lizorkin spaces is an easy consequence of this, and thus we get a new proof of the original result due to Koskela, Yang, and Zhou \cite{KYZ}.

\thm{th:f-qsinvariance} 
Let the spaces $Z$ and $Z'$ be as in the discussion above, and let $0 < s \leq 1$ and $Q/(Q+s) < q \leq \infty$. Then the composition operator induced by $\varphi$ takes $\df^{s}_{Q/s,q}(Z')$ into $\df^{s}_{Q/s,q}(Z)$, and we have
\[
  \| f\circ \varphi \|_{\df^{s}_{Q/s,q}(Z)} \lesssim \| f \|_{\df^{s}_{Q/s,q}(Z')}
\]
for all $f \in \df^{s}_{Q/s,q}(Z')$, with the implicit constant independent of $f$.
\ethm

\begin{proof}
In other words, we want to show that
\beqla{eq:mapping-condition}
  f\circ \varphi \in L^{1}_{\rm loc}(Z) \qquad \text{and} \qquad \| d\big(P(f\circ \varphi)\big) \|_{\J^{s}_{Q/s,q}(X)} \lesssim \| d(Pf) \|_{\J^{s}_{Q/s,q}(X')}
\eeq
whenever $f \in \df^{s}_{Q/s,q}(Z')$. 

For this purpose we shall estimate the quantity
\[
  v(x) := \sum_{y\sim x} \bigg(\dashint_{B(x)} |(f\circ\varphi) - (Pf\circ\Phi)(x_*)| \,d\mu + \dashint_{B(y)} |(f\circ\varphi) - (Pf\circ\Phi)(x_*)| \,d\mu \bigg)
\]
for all $x \in X$, where $x_*$ is a point at a maximal level of $X$ so that $B(x_*)$ covers $B(x)$ as well as the balls corresponding to the neighbors of $x$. Note that, a priori, $v(x)$ may be infinite for some $x$, but if $v(x)$ is finite for all $x$, then $f\circ \varphi \in L^1_{\rm loc}(Z)$ and $|d(P(f\circ \varphi))(x)| \lesssim |v(x)|$ for all $x$.

As $\varphi^{-1} : Z' \to Z$ is also quasisymmetric, the preimage of the set of Lebesgue points of $f$ has full $\mu$-measure, and for any $\xi$ in this set with points $x_{k} \in \gox \cap X_k$, $k \geq N \in \integer$, we have
\begin{align*}
  |(f \circ \varphi)(\xi) - (Pf\circ\Phi)(x_N)| \;
& \leq \;\sum_{k \geq N} |(Pf \circ \Phi)(x_{k+1}) - (Pf\circ\Phi)(x_k)| \\
& \leq \;\sum_{|y| \geq N} |d(Pf \circ \Phi)(y)| \chi_{B(y)}(\xi).
\end{align*}

In particular, if $\xi \in B(x_*)$ for some $x$, taking $x_N = x_*$ gives
\[
  v(x) \lesssim \dashint_{B(x_*)} |(f\circ\varphi) - (Pf)(\Phi(x_*))| \,d\mu \lesssim \sum_{\substack{|y| \geq |x_*| \\ B(y) \cap B(x_*) \neq \emptyset}} \frac{\mu\big(B(y)\big)}{\mu\big(B(x_*)\big)} |d(Pf \circ \Phi)(y)|.
\]
Since $|x| - \sigma \leq |x_*| \leq |x|$ for some uniform $\sigma \geq 0$, Proposition \ref{pr:basic-operator} yields
\[
  \| v \|_{\J^{s}_{Q/s,q}(X)} \lesssim \| d(Pf \circ \Phi) \|_{\J^{s}_{Q/s,q}(X)}.
\]

Now using \eqref{eq:quasiisometry}, we easily see that $|d(Pf \circ \Phi)(x)| \leq |(w\circ \Phi)(x)|$ for all $x$, where
\[
  w(x') = \sum_{y' \,:\, d_{X'}(y',x')\,\leq\,c_1} |d(Pf)(y')|
\]
for all $x' \in X'$ and $c_1 \geq 1$ is a constant. By Proposition \ref{pr:j-qsinvariance}, we thus have
$
  \| v \|_{\J^{s}_{Q/s,q}(X)} \lesssim \| w \|_{\J^{s}_{Q/s,q}(X')}.
$
Since $X'$ has bounded valency, it easily follows that
\[
  \big\|\{ 2^{|x'|s} |w(x')| : x' \in \Gamma_{1,\xi'} \} \big\|_{\ell^q} \lesssim \big\|\{ 2^{|x'|s} |d(Pf)(x')| : x' \in \Gamma_{\kappa(c_1),\xi'} \} \big\|_{\ell^q}
\]
for all $\xi' \in Z'$; so by Proposition \ref{pr:quasi-banach} (ii), we get
\[
  \| v \|_{\J^{s}_{Q/s,q}(X)} \lesssim \| d(Pf) \|_{\J^{s}_{Q/s,q}(X')}.
\]
This in particular implies that $v(x)$ is finite for all $x \in X$. As noted above, this means that $f\circ \varphi \in L^1_{\rm loc}(Z)$, and the desired estimate in \eqref{eq:mapping-condition} then follows from the inequality $|d(P(f\circ \varphi))(x)| \lesssim v(x)$ and the estimate above.
\end{proof}

We also sketch the proof of the following result, which is essentially a converse result to Theorem \ref{th:f-qsinvariance} under some mild additional assumptions on $Z$ and $Z'$. It was previously obtained for $1 < q < \infty$ in \cite[Proposition 4.3]{KKSS}. We refer to \cite{HK} for the definition of quasiconformality and the linear local connectedness property.

\prop{pr:qs-converse}
Suppose that the spaces $Z$ and $Z'$ are connected and $Q$-Ahlfors regular for some $Q > 1$, and that $Z'$ is linearly locally connected. If the composition operator induced by a homeomorphism $\psi \: Z \to Z'$ takes $\df^{s}_{Q/s,q}(Z')$ boundedly into $\df^{s}_{Q/s,q}(Z)$ for some $0 < s < 1$ and $Q/(Q+s) < q \leq \infty$, then $\psi$ is quasiconformal.
\eprop

\begin{proof}
One can follow in the argument employed in the proofs of \cite[Propositions 3.5 and 4.3]{KKSS} as long one can establish the following capacity estimates: the $\df^{s}_{Q/s,q}(Z')$-capacity of a point in $Z'$ is zero, and the $\df^{s}_{Q/s,q}(Z)$-capacity of two compact and connected subsets of $Z$ relatively close to each other is uniformly bounded from below. For Besov function spaces on connected and Ahlfors regular metric measure spaces, capacity estimates suitable for the argument mentioned above are obtained in \cite[Lemmas 3.3 and 3.4]{KKSS}. Now if we take $0 < t < s < t' < 1$, Proposition \ref{pr:jawerth-embedding} gives
\[
  \df^{t'}_{Q/t',Q/t'}(Z') \subset \df^{s}_{Q/s,q}(Z') \quad \text{and} \quad \df^{s}_{Q/s,q}(Z) \subset \df^{t}_{Q/t,Q/t}(Z)	
\]
with continuous embeddings, and since $\df^{t}_{Q/t,Q/t}(Z)$ and $\df^{t'}_{Q/t',Q/t'}(Z')$ coincide with the mentioned Besov spaces on the respective metric measure spaces with the same indices, we obtain the desired capacity estimates.
\end{proof}

Let us finally collect Theorem \ref{th:f-qsinvariance} and Proposition \ref{pr:qs-converse} in the setting of metric measure spaces where the notions of quasisymmetry and quasiconformality are known to coincide. For the fact that this indeed happens under the assumptions listed below, we refer to \cite[Theorem 5.7 and Corollary 4.10]{HK}. We again refer to \cite{HK} for the necessary definitions.

\cor{co:qs-characterization}
Suppose that $Z$ is $Q$-Ahlfors regular for some $Q > 1$, proper and quasiconvex, and that it supports a weak $(1,Q)$-Poincar\'e inequality. Suppose that $Z'$ is $Q$-Ahlfors regular, pathwise connected, locally compact and linearly locally connected. Suppose also that $Z$ and $Z'$ are both bounded or that they are both unbounded.

If $0 < s < 1$ and $Q/(Q+s) < q \leq \infty$, then the composition operator induced by a homeomorphism $\psi \: Z \to Z'$ takes $\df^{s}_{Q/s,q}(Z')$ boundedly into $\df^{s}_{Q/s,q}(Z)$ if and only if $\psi$ is quasiconformal.
\ecor

\section{Complex interpolation}\label{se:interpolation}
Finally, in the setting of Section \ref{se:unbounded} (i.e., when $Z := (Z,d,\mu)$ satisfies the doubling property \eqref{eq:doubling} and is either bounded or unbounded), we establish a complex interpolation result for our Triebel-Lizorkin spaces, which appears to be the first of its kind in the setting of doubling metric measure spaces. For this purpose we first state the result concerning the interpolation of the quasi-Banach lattices $\jp$ by the Calder\'on product method. These spaces are all embedded continuously into the space of sequences on $X$ endowed with the topology induced by the seminorms $u \mapsto |u(x)|$, $x \in X$, so any two of them are a compatible interpolation couple. Recall that in general, the Calder\'on product $X_0^{1-\theta}X_1^{\theta}$, $0 < \theta < 1$, of two quasi-Banach lattices $X_0$ and $X_1$ on a measure space $M$ is defined as the space of measurable functions $f$ on $M$ for which the quasi-norm
\[
  \|f\|_{X_0^{1-\theta}X_1^{\theta}} := \inf_{\|g\|_{X_0} \leq 1,\,\|h\|_{X_1} \leq 1} \esssup_{x \in M} \frac{|f(x)|}{|h(x)|^{1-\theta}|g(x)|^{\theta}}
\]
is finite.

The result below is analogous to \cite[Theorem 8.2]{FJ} and can be proved in a similar manner. For the reader's convenience, we have included the argument here.

\prop{pr:interpolation-j}
Let $0 < p_0 ,\, p_1 < \infty$, $0 < q_0,\, q_1 \leq \infty$, $0 < s_0,\, s_1 < \infty$, and $0 < \theta < 1$. If we define $p$, $q$, and $s$ by
\beqla{eq:interpolation-parameters}
  \frac1p = \frac{1-\theta}{p_0} + \frac{\theta}{p_1}, \quad \frac1q = \frac{1-\theta}{q_0} + \frac{\theta}{q_1} \quad \text{and} \quad s = (1-\theta)s_0 + \theta s_1,
\eeq
then
\[
  \jp = \big(\J^{s_0}_{p_0,q_0}(X)\big)^{1-\theta} \big(\J^{s_1}_{p_1,q_1}(X)\big)^{\theta}
\]
with equivalent quasi-norms.
\eprop

\begin{proof}
To show that
\[
  \big(\J^{s_0}_{p_0,q_0}(X)\big)^{1-\theta} \big(\J^{s_1}_{p_1,q_1}(X)\big)^{\theta} \subset \jp,
\]
suppose that $u \in (\J^{s_0}_{p_0,q_0}(X))^{1-\theta} (\J^{s_1}_{p_1,q_1}(X))^{\theta}$, and let $\lambda$ be a positive number strictly greater than the $(\J^{s_0}_{p_0,q_0}(X))^{1-\theta} (\J^{s_1}_{p_1,q_1}(X))^{\theta}$-quasi-norm of $u$. Thus there exist $r \in \J^{s_0}_{p_0,q_0}(X)$ and $t \in \J^{s_1}_{p_1,q_1}(X)$ so that $|u(x)| \leq \lambda |t(x)|^{1-\theta} |r(x)|^{\theta}$ for all $x$, and
\[
  \|r\|_{\J^{s_0}_{p_0,q_0}(X)} \leq 1, \quad \|t\|_{\J^{s_1}_{p_1,q_1}(X)} \leq 1.
\]
In particular,
\[
  \bigg(\sum_{x \in \gox} \big[2^{|x|s}|u(x)|\big]^q \bigg)^{1/q} \leq \lambda \bigg(\sum_{x \in \gox} \big[2^{|x|s_0}|r(x)|\big]^{(1-\theta)q} \big[2^{|x|s_1}|t(x)|\big]^{\theta q} \bigg)^{1/q}
\]
for all $\xi \in Z$ (with an obvious modification for $q = \infty$), so two applications of H\"older's inequality yield
\[
  \|u\|_{\jp} \leq \lambda \big(\|r\|_{\J^{s_0}_{p_0,q_0}(X)} \big)^{1-\theta} \big(\|t\|_{\J^{s_1}_{p_1,q_1}(X)} \big)^{\theta} \leq \lambda.
\]
Taking the infimum over admissible $\lambda$, we obtain the desired embedding.

To show that
\[
  \jp \subset \big(\J^{s_0}_{p_0,q_0}(X)\big)^{1-\theta} \big(\J^{s_1}_{p_1,q_1}(X)\big)^{\theta},
\]
we first consider the case with $q_0,\,q_1 < \infty$. Without loss of generality, assume that $p_0/q_0 \leq p_1/q_1$. Taking $u \in \jp$, we write
\[
  A_k := \bigg\{ \xi \in Z : \bigg(\sum_{x \in \gox} \big[2^{|x|s} |u(x)|\big]^q \bigg)^{1/q} > 2^k\bigg\}
\]
and
\[
  C_k := \Big\{ x \in X : \mu\big( B(x) \cap A_k \big) \geq \frac12 \mu\big(B(x)\big),\; \mu\big( B(x) \cap A_{k+1} \big) < \frac12 \mu\big(B(x)\big)\Big\}
\]
for all $k \in \integer$. Note that if $x \in X \backslash \cup_{k\in \integer} C_k$, then $u(x) = 0$. For $x \in C_k$, put
\[
  r(x) := \big( 2^{|x|v - k\gamma }|u(x)|\big)^{q/q_0} \quad {\rm and} \quad t(x) := \big( 2^{|x|w-k\delta}|u(x)|\big)^{q/q_1},
\]
where $v := s - (q_0/q)s_0$, $w := s - (q_1/q)s_1$,
\[
  \gamma := 1 - \frac{p/q}{p_0/q_0} \leq 0, \quad {\rm and} \quad \delta := 1 - \frac{p/q}{p_1/q_1} \geq 0;
\]
if $x \in X \backslash \cup_{k\in \integer} C_k$, put $r(x) = t(x) = 0$. Since $|u(x)| = |r(x)|^{1-\theta}|t(x)|^{\theta}$ for all $x$, it suffices to establish the estimates
\beqla{eq:cprod-estimate1}
  \|r\|_{\J^{s_0}_{p_0,q_0}(X)} \lesssim \|u\|_{\jp}^{p/p_0} \quad {\rm and} \quad \|t\|_{\J^{s_1}_{p_1,q_1}(X)} \lesssim \|u\|_{\jp}^{p/p_1}.
\eeq

For the first one, we note that a maximal function argument similar to the proof of Proposition \ref{pr:quasi-banach} (ii) gives
\beqla{eq:cprod-estimate2}
  \|r\|_{\J^{s_0}_{p_0,q_0}(X)}^{p_0} \approx \int_Z \bigg(\sum_{k \in \integer} \sum_{x \in \gox \cap C_k}\big[2^{|x|s_0}|r(x)|\big]^{q_0} \chi_{A_k}(\xi)\bigg)^{p_0/q_0} \,d\mu(\xi).
\eeq
If we combine this with the definition of $r$, we get
\[
  \|r\|_{\J^{s_0}_{p_0,q_0}(X)}^{p_0} \approx \int_Z \bigg(\sum_{k \in \integer} 2^{-k\gamma q}\chi_{A_k}(\xi) \sum_{x \in \gox \cap C_k}\big[2^{|x|s}|u(x)|\big]^{q} \bigg)^{p_0/q_0} \,d\mu(\xi).
\]

By the definition of $A_k$, we have $2^{k q}\chi_{A_k}(\xi) \leq \sum_{x \in \gox}\big[2^{|x|s}|u(x)|\big]^{q}$ for all $\xi \in Z$ and $k \in \integer$, and using the fact that $-\gamma \geq 0$ we obtain
\[
  \|r\|_{\J^{s_0}_{p_0,q_0}(X)}^{p_0} \lesssim \int_Z \bigg( \sum_{x \in \gox} \big[ 2^{|x|s} |u(x)|\big]^q \bigg)^{[1-\gamma](p_0/q_0)} \,d\mu(\xi) = \|u\|_{\jp}^p.
\]

The second estimate in \eqref{eq:cprod-estimate1} can be obtained in a similar manner by replacing \eqref{eq:cprod-estimate2} with
\[
  \|t\|_{\J^{s_1}_{p_1,q_1}(X)}^{p_1} \approx \int_Z \bigg(\sum_{k \in \integer} \sum_{x \in \gox \cap C_k}\big[2^{|x|s_1}|t(x)|\big]^{q_1} \chi_{(A_{k+1})^c}(\xi)\bigg)^{p_1/q_1} \,d\mu(\xi)
\]
and using the definition of $A_{k+1}$ together with the fact that $-\delta \leq 0$.

If precisely one of the numbers $q_i$ is infinite, say $q_1 < \infty = q_0$, then the argument above may be carried out by replacing $r(x)$ with $2^{(p/p_0)k - s_0|x|}$ for all $x \in C_k$. If $q_0 = q_1 = \infty$, one can take
\[
  r(x) := \big( 2^{|x|v}|u(x)|\big)^{p/p_0} \quad {\rm and} \quad t(x) := \big( 2^{|x|w}|u(x)|\big)^{p/p_1}
\]
for all $x \in X$, where $v := s - (p_0/p)s_0$ and $w := s - (p_1/p)s_1$.
\end{proof}

The Calder\'on product space $(X_0)^{1-\theta} (X_1)^{\theta}$, $0 < \theta < 1$, of two Banach lattices $X_0$ and $X_1$ is known to coincide with the corresponding interpolation space $[X_0,X_1]_\theta$ obtained by the classical complex method as long as something similar to the dominated convergence theorem holds in the resulting Calder\'on product space (i.e., if $\min(q_0,q_1) < \infty$ in the case of Proposition \ref{pr:interpolation-j} \cite[p.~125]{C}).

For quasi-Banach spaces that are \emph{A-convex} or \emph{analytically convex}, an extension of the classical complex interpolation theory of Banach spaces has been developed by Kalton et al., see for instance \cite[Section 7]{KMM} and the references therein. In our case, the quasi-norm of $\jp$ satisfies
\[
  \Big\| \Big( \sum_{1 \leq i \leq n} |u_i|^r\Big)^{1/r}\Big\|_{\jp} \leq \Big( \sum_{1 \leq i \leq n} \|u_i\|_{\jp}^r \Big)^{1/r}
\]
for all sequences $u_1,\dots, u_n$ on $X$ and $0 < r \leq \min(1,p,q)$; so $\jp$ is $A$-convex for all $p$, $q$, and $s$, as is the sum of any two such spaces. This in particular means that a Calder\'on product of two such spaces will coincide with the corresponding complex interpolation space. We refer, for example, to \cite[Theorem 4.4]{K}, the remark following \cite[Proposition 10]{MM}, and \cite[Section 7]{KMM} and the references therein for details.

We can now state the interpolation result for our Triebel-Lizorkin spaces. 

\thm{th:interpolation-f}
Suppose that $0 < s_0,\;s_1 < 1$, $Q/(Q+s_0) < p_0,\, q_0 < \infty$, $Q/(Q+s_1) < p_1 < \infty$, and $Q/(Q+s_1) < q_1 \leq \infty$. For $0 < \theta < 1$ we then have
\[
  \big[\df^{s_0}_{p_0,q_0}(Z),\df^{s_1}_{p_1,q_1}(Z)\big]_{\theta} = \dfp
\]
with $p$, $q$, and $s$ as in \eqref{eq:interpolation-parameters}.
\ethm

A similar result with $p_i$, $q_i \in (1,\infty)$ in the setting of Ahlfors regular spaces can be found in \cite{HS}.    

The proof of Theorem \ref{th:interpolation-f} is a standard application of Proposition \ref{pr:interpolation-j} above coupled with the method of retractions and co-retractions \cite[Lemma 7.11]{KMM}. In our situation the construction of a suitable retraction from a sequence space of the type $\J^s_{p,q}$ to a function space of the type $\df^s_{p,q}$ is a nontrivial task, and because the existence of such a retraction turns out to be very useful  in other contexts as well (see \cite{SS}), we will state it as a separate result. The proof of Theorem \ref{th:interpolation-f} will be postponed after the proof of Proposition \ref{pr:retraction} below.

In what follows it will be convenient to work with sequence spaces on the edges $E$, as the discrete derivative of a function on $X$ is naturally defined as a scalar-valued function on $E$. We equip the edges in $E$ with an orientation and denote by $e_{x,x'}$ the directed edge from $x$ to $x'$ for any neighbors $x$, $x' \in X$. The orientation is chosen in such a way that if $x\sim x'$ and $|x'| > |x|$, then $x'$ is the endpoint of the edge joining $x$ and $x'$. For a sequence $v\colon X\to\complex$, we write $dv$ for the sequence defined on $E$ such that $dv(e_{x,x'}) = v(x') - v(x)$ for all admissible $x$, $x' \in X$.

To give a concrete definition of the spaces $\jpe$ of sequences defined on $E$, we need some additional notation. For an edge $e \in E$ joining the vertices $x \in X$ and $x' \in X$, let $|e| := |x|\land |x'|$ and $B(e) := B(x) \cup B(x')$. For $0 < p < \infty$, $0 < q \leq \infty$ and $0 < s < \infty$, write $\jpe$ for the space of sequences $u$ on $E$ for which
\[
  \|u\|_{\jpe} := \Big( \int_{Z} \big\| \{ 2^{|e|s} |u(e)|\chi_{B(e)}(\xi)\}\big\|_{\ell^q(E)}^p \,d\mu(\xi)\Big)^{1/p}
\]
is finite. Proposition \ref{pr:interpolation-j} obviously continues to hold with $E$ in place of $X$.

\prop{pr:retraction}
Suppose that $0 < s < 1$, $Q/(Q+s) < p < \infty$ and $Q/(Q+s) < q \leq \infty$. Then there exist bounded linear operators
\[
  S\colon \dfp\to\jpe \quad \text{and} \quad R\colon \jpe\to \dfp
\]
such that $R\circ S$ is the identity mapping on $\dfp$. More explicitly, if  $\xi_0$ is an arbitrary fixed point of $Z$, we may take
\[
  S := f \mapsto d(Pf)
\]
and
\[
  R := u \mapsto \lim_{N\to\infty} \bigg( \sum_{n=-N}^{N} I_n u(\cdot) - \sum_{n=-N}^{-1} I_n u(\xi_0)\bigg) + \complex,
\]
where
\[
  I_n u := \sum_{(y,y') \in (X_n \times X_{n+1}), \; y\sim y' }u(e_{y,y'})\psi_y \psi_{y'}
\]
for all $n \in \integer$.
\eprop

\begin{proof}
The boundedness of $S$ is clear (from the definition of the spaces $\dfp$), so our main task is to verify that the operator $R$ is well-defined and bounded between the desired spaces, and that $R\circ S$ is the identity operator on $\dfp$. In a sense we need to ``integrate'' arbitrary sequences in $\jpe$ which do not necessarily come from discrete derivatives of sequences defined on $X$.

For simplicity, let us first consider the case where $Z$ is bounded. Without loss of generality, we may assume that $\diam Z = 1$ as in Section \ref{se:definitions}, and we are only concerned with the vertices of $X$ with $|x| \geq 0$. Then 
\[
  \integral u := \sum_{n = 0}^{\infty} I_n u
\]
converges in $L^1_{\rm loc}(Z)$ and pointwise $\mu$-almost everywhere, as can be shown by an argument similar to the proof of Lemma \ref{le:trace}. Since we have

\[
  \sum_{n=0}^{N-1} I_n u(\xi) = \sum_{ \substack{(y_0,\cdots,y_N) \in (X_0 \times \cdots \times X_N) \\ y_0\sim y_1\sim \cdots \sim y_N}} \bigg( \sum_{i=1}^N u(e_{y_{i-1},y_i})\bigg) \psi_{y_0}(\xi)\cdots\psi_{y_N}(\xi)
\]
for all $N \in \nanu$, the expression $\integral u(\xi)$ can be thought of as a weighted average of the ``integrals'' of $u$ along the geodesic segments (starting at $\bzero$) of the metric graph $(X,E)$ contained in a cone at $\xi$. Also, if $v$ is a sequence on $X$ such that $|dv|$ is in $\jp$, we obviously have $I_n(dv) = T_{n+1}v - T_n v$, and thus $\integral (dv) = \trace v - v(\bzero)$. So taking
\[
  R := u \mapsto \integral u + \complex,
\]
$R\circ S$ is the identity mapping on $\dfp$. We thus have to check that $R$ takes $\jpe$ continuously into $\dfp$.

Let $u \in \jpe$ and write
\[
  U_k := \bigg( \sum_{|e| = k} \big[2^{|e|s} |u(e)|\big]^{q} \chi_{B(e)} (\cdot)\bigg)^{1/q}
\]
for all $k \in \nanu_0$, so that the quasi-norm of $u$ is obtained as the $L^{p}(\ell^{q})$-quasi-norm of the functions $U_k$.

By the $L^1_{\rm loc}$-convergence, we have
\[
  2^{|x|s}|d(P(\integral u))(x)| \leq \sum_{n \geq 0} 2^{|x|s}|d(P(I_n u))(x)|
\]
for all $x \in X$. Suppose $x \in \gox$ for some fixed $\xi \in Z$. We first estimate the terms of the series above with $n \leq |x|$. If $x'$ is any neighbor of $x$, we have
\begin{align*}
& 2^{|x|s} |P(I_n u)(x) - P(I_n u)(x')| \\
\leq \; & 2^{|x|s} \sum_{\substack{|e_{y,y'}| = n \\ B(e_{y,y'})\cap (B(x)\cup B(x')) \neq \emptyset}} |u(e_{y,y'})|\dashint_{B(x)}\dashint_{B(x')} |\psi_y (\eta)\psi_{y'}(\eta) - \psi_y(\eta') \psi_{y'}(\eta')| \,d\mu(\eta) \,d\mu(\eta').
\end{align*}

Using the Lipschitz continuity of the functions $\psi_y$ and $\psi_{y'}$ above, as well as the fact that $$\#\{ e : |e|=n\;\textrm{and}\; B(e) \cap (B(x)\cup B(x')) \neq \emptyset\}$$ is bounded uniformly in $n \leq |x|$, we get 
\begin{align*}
  2^{|x|s} |P(I_n u)(x) - P(I_n u)(x')| \;
& \lesssim \;2^{(n - |x|)(1 - s)} \sum_{\substack{|e| = n \\ B(e) \cap (B(x)\cup B(x')) \neq \emptyset}} 2^{|e|s}|u(e)| \\
& \lesssim \;2^{(n - |x|)(1 - s)} \hlmax\big(U_n^r\big)(\xi)^{1/r},
\end{align*}
where $r$ is chosen so that $Q/(Q+s) < r < \min(1,p,q)$. As $X$ has bounded valency, we infer that $2^{|x|s}|d(P(I_n u))(x)|$ is bounded from above by a constant times the rightmost quantity above.

Now if $n > |x|$, let $x_*$ be a point at a maximal level of $X$ so that $B(x_*)$ covers the balls corresponding to $x$ and its neighbors. Then $|x| - \sigma \leq |x_*| \leq |x|$ for some uniform constant $\sigma \geq 0$, and we have
\begin{align*}
  2^{|x|s}|d(P(I_n u))(x)| \;
& \lesssim \;2^{(|x|-n)s} \dashint_{B(x_*)} \sum_{|e| = n} 2^{|e|s} |u(e)|\chi_{B(e)}(\eta) \,d\mu(\eta) \\
& \lesssim \;2^{(|x|-n)s} \sum_{\substack{|e| = n \\ B(e) \cap B(x_*) \neq \emptyset}} \frac{\mu(B(e))}{\mu(B(x_*))} 2^{|e|s} |u(e)| \\
& \lesssim \;2^{(|x|-n)(Q+s - Q/r)} \hlmax\big( U_n^r)(\xi)^{1/r};
\end{align*}
see \eqref{eq:sk-1} and \eqref{eq:sk-2} in the proof of Proposition \ref{pr:basic-operator}.

Combining these estimates and writing $\lambda := \min(Q+s-Q/r , 1-s)(q \land 1) > 0$, we get
\beqla{eq:integration-pointwise}
  \big[2^{|x|s}|d(P(\integral u))(x)|\big]^{q} \lesssim \sum_{n \geq 0} 2^{-\lambda |n-|x||} \hlmax\big( U_n^r)(\xi)^{q/r},
\eeq
and summing over $x \in \gox$, we further get
\[
  \sum_{x \in \gox} \big[2^{|x|s}|d(P(\integral u))(x)|\big]^{q} \lesssim \sum_{n \geq 0} \hlmax\big( U_n^r)(\xi)^{q/r}.
\]
One can then finish using the Fefferman-Stein maximal theorem.

Let us now consider the case where $Z$ is either bounded or unbounded, and $X_k$ is defined for all $k \in \integer$. We fix a point $\xi_0 \in Z$. Now, for $u \in \jpe$, put
\[
  \integral u(\xi) := \lim_{N\to\infty} \bigg(\sum_{n=-N}^{N}I_nu(\xi) - \sum_{n = -N}^{-1}I_n u(\xi_0)\bigg).
\]
To justify the existence of this limit (in $L^1_{\rm{loc}}(Z)$ and pointwise almost everywhere), we first observe that an argument similar to the proof of Lemma \ref{le:trace} again yields the convergence of $\sum_{n\geq 0} I_n u$ in $L^1_{\rm{loc}}(Z)$. In order to treat the remaining part of the sum, we note that for a fixed $k \in \integer$ and all integers $n \leq k$, the Lipschitz continuity of the functions $\psi_y$ yields
\[
  \int_{B(\xi_0,2^{-k})} |I_n u(\xi) - I_n u(\xi_0)| \,d\mu(\xi) \lesssim \mu\big(B(\xi_0,2^{-k})\big)\sum_{\substack{|e| = n \\ B(e) \cap B(\xi_0,2^{-k})\neq\emptyset}} 2^{n-k}|u(e)|.
\]
Writing $4B(e)$ for $4B(x) \cup 4B(x')$ for any edge $e$ joining the points $x$, $x' \in X$, we thus have

\begin{align*}
& \sum_{n \leq k} \int_{B(\xi_0,2^{-k})} |I_n u(\xi) - I_n u(\xi_0)| \,d\mu(\xi) \\
& \qquad \lesssim \;c(\xi_0,k) \sum_{\substack{|e| \leq k \\ B(e)\cap B(\xi_0,2^{-k})\neq\emptyset}} 2^{|e|} |u(e)| \\
& \qquad \lesssim \;c(\xi_0,k) \inf_{\xi \in B(\xi_0,2^{-k})}\sum_{|e| \leq k} 2^{|e|} |u(e)|\chi_{4 B(e)}(\xi) \\
& \qquad \lesssim \;c(\xi_0,k,q,s) \inf_{\xi \in B(\xi_0,2^{-k})} \big\|\{2^{|e|s} |u(e)| \chi_{4B(e)}(\xi): |e| \leq k\} \|_{\ell^{q}}.
\end{align*}
Here the last quantity is finite, since $\xi \mapsto \big\|\{2^{|e|s} |u(e)| \chi_{4B(e)}(\xi): |e| \leq k\} \|_{\ell^{q}}$ is a function in $L^{p}(Z)$, as can be seen by an argument similar to the proof of Proposition \ref{pr:quasi-banach} (ii).

In particular, $\integral (dv) = \trace v - T_0 v(\xi_0)$ almost everywhere for all $v$ such that $|dv|$ is in $\jp$. So taking $R$ as in the statement of the result, $R\circ S$ is the identity operator on $\dfp$, and $R$ is bounded between the desired spaces, as can be checked by the same argument as in the case where $Z$ is bounded, with \eqref{eq:integration-pointwise} replaced by
\[
  \big[2^{|x|s}|d(P(\integral u))(x)|\big]^{q} \lesssim \sum_{n \in \integer} 2^{-\lambda |n-|x||} \hlmax\big( U_n^r)(\xi)^{q/r} \qedhere
\]
\end{proof}

We are now ready to give the proof of our main interpolation theorem.

\begin{proof}[Proof of Theorem \ref{th:interpolation-f}]
Write $U$ for the quasi-Banach space $\J^{s_0}_{p_0,q_0}(E)$ + $\J^{s_1}_{p_1,q_1}(E)$ of sequences on $E$. More precisely, $U$ is defined as the collection of sequences on $E$ that can be expressed as $u_0 + u_1$ with $u_i \in \J^{s_i}_{p_i,q_i}(E)$. The quasi-norm of $u$ is defined as the infimum of $\|u_0\|_{\J^{s_0}_{p_0,q_0}(E)} + \|u_1\|_{\J^{s_1}_{p_1,q_1}(E)}$ over admissible representations. The spaces $\J^{s_i}_{p_i,q_i}(E)$ are then both continuously embedded into $U$ in a natural way. Similarly, we write $V$ for the quasi-Banach space $\df^{s_0}_{p_0,q_0}(Z) + \df^{s_1}_{p_1,q_1}(Z)$ of (equivalence classes of) functions in $L^1_{\rm loc}(Z) / \complex$.

It is easily seen that the operators $S$ and $R$ in Proposition \ref{pr:retraction} extend as operators $S \: V \to U$ and $R \: U \to V$ so that $R\circ S$ is the identity mapping on $V$, $S$ takes $\df^{s_i}_{p_i,q_i}(Z)$ continuously into $\J^{s_i}_{p_i,q_i}(E)$ for $i \in \{0,1\}$ and $R$ takes $\J^{s_i}_{p_i,q_i}(E)$ continuously into $\df^{s_i}_{p_i,q_i}(Z)$ for $i \in \{0,1\}$. Using the method of retractions and co-retractions \cite[Lemma 7.11]{KMM}, we can thus deduce that
\[
  \big[\df^{s_0}_{p_0,q_0}(Z),\df^{s_1}_{p_1,q_1}(Z)\big]_{\theta} = R\big(\big[\J^{s_0}_{p_0,q_0}(Z),\J^{s_1}_{p_1,q_1}(Z)\big]_{\theta}\big) = R\big(\jp\big) = \dfp
\]
with equivalent quasi-norms.
\end{proof}

Let us finally return to the proof of Corollary \ref{co:lipschitz-density} in the case that $Z$ is unbounded, which was postponed until now.

\begin{proof}[Completion of the proof of Corollary \ref{co:lipschitz-density}]
Suppose that $Z$ is unbounded, $0 < s < 1$ and $Q/(Q+s) < p,q < \infty$. The set $F$ of sequences on $E$ with finite support is obviously dense in $\jpe$. Taking $R \: \jpe \to \dfp$ as in Proposition \ref{pr:retraction}, we thus have that $R(F)$ is a dense subspace of $\dfp$. By construction, the elements of $R(F)$ are obviously (equivalence classes of) Lipschitz functions with bounded support.
\end{proof}

\end{document}